\documentclass[12pt]{article}
\usepackage{amssymb}

\newcommand{\real}{{\bf R}}
\newcommand{\complex}{{\bf C}}
\newcommand{\intplus}{{\bf N}}

\newcommand{\e}{{\rm e}}                
\renewcommand{\d}{\,{\rm d}}            
\renewcommand{\Re}{{\rm Re\,}}

\newcommand{\loc}{{\rm loc}}

\renewcommand{\epsilon}{\varepsilon}
\renewcommand{\phi}{\varphi}

\newcommand{\cB}{{\cal B}}

\newcommand{\cL}{{\cal L}}

\newcommand{\cM}{{\cal M}}
\newcommand{\cO}{{\cal O}}

\newcommand{\cG}{{\cal G}}
\newcommand{\cT}{{\cal T}}

\newcommand{\uu}{{\bf u}}
\newcommand{\xx}{{\bf x}}
\newcommand{\yy}{{\bf y}}
\newcommand{\VV}{{\bf V}}
\newcommand{\oone}{{\bf 1}}
\newcommand{\LL}{{\mathrm L}}
\newcommand{\UU}{{\bf U}}
\newcommand{\OO}{{\bf \Omega}}

\newcommand{\OB}{{\Omega^B}}
\newcommand{\oomega}{{\mbox{\boldmath$ \omega$}}}

\newtheorem{theorem}{Theorem}[section]

\newtheorem{proposition}[theorem]{Proposition}
\newtheorem{corollary}[theorem]{Corollary}
\newtheorem{remark}[theorem]{Remark}

\newcommand{\PPP}{{\mathbb P}}
\newcommand{\LLL}{{\mathbb L}}
\newcommand{\NNN}{{\mathbb N}}
\newcommand{\HHH}{{\mathbb H}}
\newcommand{\XXX}{{\mathbb X}}
\newcommand{\YYY}{{\mathbb Y}}
\newcommand{\FFF}{{\mathbb F}}

\newcommand{\reff}[1]{(\ref{#1})}

\newcommand{\inttwo}{\int_{\real^2}}
\newcommand{\proof}{{\noindent \bf Proof:\ }}

\newcommand{\half}{{\frac{1}{2}}}
\newcommand{\xp}{{x_{\perp}}}
\newcommand{\yp}{{y_{\perp}}}
\newcommand{\cLp}{{\cL}_{\perp}}

\def\build#1_#2^#3{\mathrel{
  \mathop{\kern 0pt#1}\limits_{#2}^{#3}}}
\def\QED{\mbox{}\hfill$\Box$}

\usepackage{amsfonts}
\begin{document}

\title{Three-Dimensional Stability of Burgers Vortices: the Low
Reynolds Number Case.}

\author{Thierry Gallay \\ Institut Fourier \\
Universit\'e de Grenoble I \\ BP 74 \\
38402 Saint-Martin-d'H\`eres \\France
\and
C. Eugene Wayne \\
Department of Mathematics \\
and  Center for BioDynamics \\
Boston University \\
111 Cummington St.\\
Boston, MA 02215, USA}

\date{\normalsize March 15, 2005}

\maketitle
\begin{abstract} 
In this paper we establish rigorously that the family of Burgers
vortices of the three-dimensional Navier-Stokes equation is stable for
small Reynolds numbers. More precisely, we prove that any solution
whose initial condition is a small perturbation of a Burgers vortex
will converge toward another Burgers vortex as time goes to infinity,
and we give an explicit formula for computing the change in the
circulation number (which characterizes the limiting vortex completely.)  
We also give a rigorous proof of the existence and stability of
non-axisymmetric Burgers vortices provided the Reynolds number is
sufficiently small, depending on the asymmetry parameter.
\end{abstract}

\section{Introduction}\label{intro} 

Numerical simulations of turbulent flows have lead to the general
conclusion that vortex tubes serve as important organizing structures
for such flows -- in the memorable phrase of \cite{moffatt:1994} they
form the ``sinews of turbulence''.  After the discovery by Burgers
\cite{burgers:1948} of the explicit vortex solutions of the
three-dimensional Navier-Stokes equation which now bear his name,
these solutions have been used to model various aspects of turbulent
flows \cite{townsend:1951}.  It was also observed in numerical
computations of fluid flows that the vortex tubes present in these
simulations usually did not exhibit the axial symmetry of the explicit
Burgers solution, but rather an elliptical core region. This lead to a
search for non-axisymmetric vortices \cite{robinson:1984},
\cite{moffatt:1994}, \cite{jimenez:1996}. While no rigorous proof of
their existence was available until recently, perturbative
calculations and extensive numerical simulations have lead to the
expectation that stationary vortical solutions of the
three-dimensional Navier-Stokes equation do exist for any Reynolds
number and all values of the asymmetry parameter (which we define
below) between zero and one.

When addressing the stability of Burgers vortices, it is very
important to specify the class of allowed perturbations. If we
consider just two-dimensional perturbations (i.e., perturbations which
do not depend on the axial variable), then fairly complete answers are
known. Robinson and Saffman \cite{robinson:1984} computed
perturbatively the eigenvalues of the linearized operator at the
Burgers vortex and proved its stability for sufficiently small
Reynolds numbers. Numerical computations of these eigenvalues were
performed by Prochazka and Pullin \cite{prochazka:1995}, and no
instability was found up to $Re = 10^4$. A similar conclusion was
drawn for non-symmetric vortices \cite{prochazka:1998}. The first
mathematical work is \cite{gallay:2004}, where we proved that the
axisymmetric Burgers vortex is {\it globally stable} with respect to
integrable, two-dimensional perturbations, for any value of the
Reynolds number.  Decay rates in time of spatially localized
perturbations were also computed, explaining partially the numerical
results of \cite{prochazka:1995}. Building on this work the existence
and local stability of slightly asymmetric vortices with respect to
two dimensional perturbations was proved in \cite{gallay:2005} for
arbitrary Reynolds numbers.

The stability issue is much more difficult if we allow for
perturbations which depend on the axial variable too, and very few
results have been obtained so far in this truly three-dimensional
case. One early study by Leibovich and Holmes \cite{leibovich:1981}
concluded that one could not prove global stability for any Reynolds
number solely by means of energy methods.  Using a kind of Fourier
expansion in the axial variable, Rossi and Le Diz\`es
\cite{rossi:1997} showed that the point spectrum of the linearized
operator is associated with purely two-dimensional perturbations.
Crowdy \cite{crowdy:1998} obtained a formal asymptotic expansion of
the eigenfunctions in the axial variable. In an important recent work,
Schmid and Rossi \cite{schmid:2004} rewrote the linearized equations
in a form which allowed them to compute numerically the evolution of
various Fourier modes, from which they concluded that eventually all
perturbative modes will be damped out.

In this paper we address rigorously the existence of
non-axisymmetric vortices and the stability with respect to 
three-dimensional perturbations of both the symmetric and non-symmetric
vortex solutions. More precisely we will prove that, for all
values of the asymmetry parameter between zero and one, 
non-axisymmetric vortices exist at least for small Reynolds numbers.
In addition, we show that this family of vortex solutions is, in 
the language of dynamical systems theory, {\em asymptotically stable
with shift}. That is to say, if we take initial conditions that
are small perturbations of a vortex solution, the resulting solution
of the Navier-Stokes equation will converge toward a vortex solution, 
but not, in general, the one which we initially perturbed. We also 
give a formula for computing the limiting vortex toward which the 
solution converges.

\medskip
We now state our results more precisely. The three-dimensional 
Navier-Stokes equations for an incompressible fluid with constant 
density $\bar{\rho}$ and kinematic viscosity $\nu$ are the partial 
differential equations:
\begin{equation}\label{3DNS}
  \partial_t \uu + (\uu \cdot \nabla) \uu \,=\, \nu \Delta \uu - 
  \frac{1}{\bar{\rho}}\nabla p~, \quad \nabla \cdot \uu \,=\, 0~.
\end{equation}
Here $\uu(x,t)$ is the velocity of the fluid and $p(x,t)$ its
pressure. Equation~\reff{3DNS} will be considered in the whole 
space $\real^3$. Burgers vortices are particular solutions of 
\reff{3DNS} which are perturbations of the background straining 
flow
\begin{equation}\label{strain}
  \uu^s(x) \,=\, \pmatrix{\gamma_1 x_1 \cr \gamma_2 x_2 \cr
  \gamma_3 x_3}~, \quad p^s(x) \,=\, -\frac12 \bar\rho 
  (\gamma_1^2 x_1^2 + \gamma_2^2 x_2^2 + \gamma_3^2 x_3^2)~,
\end{equation}
where $\gamma_1, \gamma_2, \gamma_3$ are real constants satisfying
$\gamma_1 + \gamma_2 + \gamma_3 = 0$. We restrict ourselves to 
the case of an  {\em axial strain} aligned with the vertical axis, 
namely we assume $\gamma_1, \gamma_2 < 0$ and $\gamma_3 > 0$. 
Setting $\uu = \uu^s + \UU$, we obtain the following evolution 
equation for the vorticity $\OO = \nabla \times \UU$:
\begin{equation}\label{3DV1}
  \partial_t \OO + (\UU\cdot \nabla) \OO - (\OO\cdot\nabla)\UU
   + (\uu^s \cdot \nabla) \OO - (\OO\cdot\nabla)\uu^s
  \,=\, \nu\Delta\OO~, \quad \nabla \cdot \OO \,=\, 0~.
\end{equation}
Under reasonable assumptions which will be satisfied for the solutions
we consider, the rotational part $\UU$ of the velocity can be recovered
from the vorticity $\OO$ by means of the Biot-Savart law:
\begin{equation} \label{BS3}
  \UU(x) = -\frac{1}{4\pi} \int_{\real^3} 
  \frac{(\xx-\yy) \times \OO(y)}{|x - y|^3} \d y\ ,
  \quad x \in \real^3~.
\end{equation}

For notational simplicity we begin by discussing the {\em axisymmetric}
case where $\gamma_1 = \gamma_2 = - \gamma_3/2$. In this situation, 
it is well-known \cite{burgers:1948} that \reff{3DV1} has a family of 
explicit stationary solutions of the form $\OO = \Gamma \hat \OO^B$, 
where $\Gamma \in \real$ is a parameter and 
\begin{equation}\label{hatOBdef}
  \hat \OO^B(\xp) \,=\, \pmatrix{0 \cr 0 \cr \hat\Omega^B(\xp)}~, \quad 
  \hat\Omega^B(\xp) \,=\, \frac{\gamma }{4 \pi \nu}
  \,e^{-\gamma |\xp|^2/(4\nu)}~.
\end{equation}
Here $\xp = (x_1,x_2)$, $|\xp|^2 = x_1^2 + x_2^2$, and $\gamma \equiv
\gamma_3 > 0$. The velocity field corresponding to $\Gamma\hat\OO^B$ 
is $\Gamma\hat\UU^B$, where
\begin{equation}\label{hatUUBdef}
  \hat \UU^B(\xp) \,=\, \frac{1}{2\pi} \pmatrix{-x_2 \cr x_1 \cr 0}
  \frac{1}{|\xp|^2} \Bigl(1 - e^{-\gamma |\xp|^2/(4\nu)} \Bigr)~.
\end{equation} 
These solutions are called the {\em axisymmetric Burgers vortices}. 
Observe that $\inttwo \hat\Omega^B\d\xp = 1$, so that the parameter 
$\Gamma$ represents the circulation of the velocity field $\Gamma 
\hat \UU^B$ at infinity (in the horizontal plane $x_3 = 0$). 
Following \cite{moffatt:1994}, we define the Reynolds number
associated with the Burgers vortex $\Gamma \hat \UU^B$ as
\begin{equation}\label{Rdef}
  R \,=\, \frac{|\Gamma|}{\nu}~.
\end{equation}

Our principal result concerns the evolution of solutions of
\reff{3DV1} with initial conditions that are close to a Burgers
vortex.  Unlike in much previous work the perturbations we consider do
not merely depend on the transverse variables $\xp$, but also on
$x_3$. We prove that any solution of \reff{3DV1} starting sufficiently
close to the Burgers vortex with circulation $\Gamma$ converges as 
$t \to +\infty$ toward a Burgers vortex with circulation $\Gamma'$ 
close to $\Gamma$, and we give an explicit formula for computing
the difference $\Gamma' - \Gamma$ in terms of the initial 
perturbation.

We now introduce some function spaces to measure the size of our 
perturbations. Roughly speaking, we require the perturbations
to decay as inverse powers of $\xp$ as $|\xp| \to \infty$, 
but need only boundedness in $x_3$. To be specific, we use 
the following weight function
\begin{equation}\label{rhodef}
  b(\xp) \,=\, (1+\gamma x_1^2/\nu)^{1/2} 
  (1+\gamma x_2^2/\nu)^{1/2}~, \quad \xp = (x_1,x_2) \in \real^2~. 
\end{equation}
Given any $m \ge 0$, we define $L^2(m) \,=\, \{\omega : \real^2 \to \real 
\,|\, \|\omega\|_{L^2(m)} < \infty\}$, where 
\begin{equation}\label{L2mdef}
  \|\omega\|_{L^2(m)}^2 \,=\, \frac{1}{\gamma \nu}\inttwo b(\xp)^{2m} 
  |\omega(\xp)|^2 \d\xp~.
\end{equation}
In other words, a function $\omega$ belongs to $L^2(m)$ if and only
if $\omega$, $|x_1|^m \omega$, $|x_2|^m \omega$, and $|x_1 x_2|^m 
\omega$ are square integrable over $\real^2$. For later use, we
observe that $L^2(m)$ is continuously embedded into $L^1(\real^2)$
if $m > 1/2$, i.e. there exists $C > 0$ such that $\|\omega\|_{L^1} 
\le C \|\omega\|_{L^2(m)}$ for all $\omega \in L^2(m)$. 

Our main space $X^2(m)$ will be the set of all $\omega : \real^3 
\to \real$ such that $\xp \mapsto \omega(\xp,x_3) \in L^2(m)$ 
for all $x_3 \in \real$, and such that the map $x_3 \mapsto 
\omega(\cdot,x_3)$ is bounded and continuous from $\real$ into
$L^2(m)$. As is easily verified, $X^2(m) \simeq C^0_b(\real,L^2(m))$
is a Banach space equipped with the norm
\begin{equation}\label{X2mdef}
  \|\omega\|_{X^2(m)} \,=\, \sup_{x_3 \in \real}
  \|\omega(\cdot,x_3)\|_{L^2(m)}~.
\end{equation}

\begin{remark} If $\oomega = (\omega_1,\omega_2,\omega_3)$ is a vector
field whose components are elements of $X^2(m)$, we shall often
write $\oomega \in X^2(m)$ instead of $\oomega \in X^2(m)^3$, and
$\|\oomega\|_{X^2(m)}$ instead of $\|(\omega_1^2 + \omega_2^2
+\omega_3^2)^{1/2}\|_{X^2(m)}$. A similar abuse of notation will 
occur for other function spaces too.
\end{remark}

Consider initial conditions for the vorticity equation \reff{3DV1}
which are a perturbation of the Burgers vortex:
\begin{equation}\label{IC1}
  \OO^0(x) \,=\, \Gamma \hat \OO^B(\xp) + \oomega^0(x)~,
\end{equation}
with $\Gamma \in \real$ and $\oomega^0 \in X^2(m)^3$. Define
\begin{equation}\label{dRdef}
  \delta \Gamma \,=\, \Bigl(\frac{\gamma}{2\pi\nu}\Bigr)^{1/2}
  \int_{\real} \inttwo e^{-\gamma x_3^2/(2\nu )}
  \omega^0_3(\xp,x_3) \d\xp \d x_3~.
\end{equation}
Just to make sure the notation is clear, in the integrand
$\omega^0_3$ refers to the third component of the initial perturbation
$\oomega^0$.

\begin{theorem}\label{sym_stable} Fix $m > 3/2$, and assume that
$(\gamma_1,\gamma_2,\gamma_3) = \gamma(-\half,-\half,1)$. For 
any $\mu \in (0,1/2)$, there exist $R_0 > 0$ and $\epsilon_0 > 0$
such that if $|\Gamma| \le R_0 \nu$ and $\|\oomega_0\|_{X^2(m)} \le 
\epsilon_0$, then the solution $\OO(x,t)$ of \reff{3DV1} with initial 
condition \reff{IC1} converges as $t$ tends to infinity to the Burgers
vortex with circulation number $\Gamma'= \Gamma + \delta \Gamma$,
where $\delta \Gamma$ is given by \reff{dRdef}.  Convergence is with
respect to the $L^2(m)$ norm in $\xp$ and uniformly on compact sets in
$x_3$.  More explicitly, if $I \subset \real$ is any compact interval
we have
\begin{equation}\label{conv1}
  \sup_{x_3 \in I} \|\OO(\cdot,x_3,t) - \Gamma'\hat\OO^B
  (\cdot)\|_{L^2(m)} \,=\, \cO(e^{-\mu \gamma t})~, \quad t \to +\infty~.
\end{equation}
\end{theorem}

\begin{remark}
Here and in the sequel, all constants are independent of the 
strain $\gamma$ and the viscosity $\nu$. In fact, both parameters 
will shortly be eliminated by a rescaling.
\end{remark}

The proof of Theorem \ref{sym_stable} uses ideas from our analysis of
the stability of the two-dimensional Oseen vortices in
\cite{gallay:2004}. The main observation is that, if we linearize
equation \reff{3DV1} at the Burgers vortex $\Gamma\hat\OO^B$ for
small $\Gamma$, we obtain a small perturbation of a non-constant
coefficient differential operator for which we can explicitly compute
an integral representation of the associated semigroup. This semigroup
decays exponentially when acting on functions $\oomega \in X^2(m)^3$
{\em provided} $\omega_3 \in X^2_0(m)$, where
\begin{equation}\label{X20mdef}
  X^2_0(m) \,=\, \Bigl\{\omega \in X^2(m) ~\Big|~ \inttwo \omega(\xp,x_3) 
  \d\xp = 0 \quad {\rm for~all~} x_3 \in \real\Bigr\}~.
\end{equation}
Thus an important step in the proof consists in decomposing the 
perturbation as $\oomega(x,t) = \phi(x_3,t)\hat\OO^B(\xp) + 
\tilde \oomega(x,t)$, where
$$
  \phi(x_3,t) \,=\, \inttwo \omega_3(\xp,x_3,t)\d\xp~, 
  \quad x_3 \in \real~, \quad t \ge 0~.
$$ 
By construction $\tilde \omega_3 \in X^2_0(m)$, hence $\tilde 
\oomega(x,t)$ will decay exponentially to zero by the remark above.
Now, the crucial point is that $\phi(x_3,t)$ satisfies the
amazingly simple equation
$$
  \partial_t \phi + \gamma x_3\partial_3 \phi \,=\, \nu\partial_3^2
  \phi~,
$$
which can be solved explicitly, see \reff{phi_solve} below. From the
solution formula we see that $\phi(x_3,t)$ converges uniformly on
compact sets to the value $\delta\Gamma$ as $t \to \infty$, and
\reff{conv1} follows.  In other words, $\phi$ is a ``zero mode'' which 
is responsible for the fact that the family of Burgers vortices is only
asymptotically stable {\em with shift}.

\medskip
While the axisymmetric vortex solution has been extensively studied because
of the explicit formulas for its velocity and vorticity fields, 
numerical experiments on turbulent flows seem to indicate that the vortex 
tubes that are prominent in these flows are not symmetric, but rather 
elliptical in cross section. A natural way to obtain such vortices
is to assume that the straining flow is not axisymmetric, namely
\begin{equation}\label{nonsymgam}
  \gamma_1 \,=\, -\frac{\gamma}2 (1+\lambda)~, \quad 
  \gamma_2 \,=\, -\frac{\gamma}2 (1-\lambda)~, \quad 
  \gamma_3 \,=\, \gamma~, 
\end{equation}
where $\gamma > 0$ and $\lambda \in [0,1)$ is an additional parameter
which measures the asymmetry of the strain. While no explicit formulas
for the vortex are known when $\lambda > 0$, extensive perturbative and
numerical investigations indicate that there do exist stationary
solutions of \reff{3DV1} for $0 \le \lambda < 1$, which for $\lambda$
close to zero are small perturbations of the vorticity field of the
axisymmetric Burgers vortex \cite{robinson:1984}, \cite{moffatt:1994},
\cite{prochazka:1998}.  As in the symmetric case, there is in fact a
family of vortices for each value of $\lambda \in [0,1)$ parametrized
by the total circulation $\Gamma$, but when $\lambda > 0$ these
solutions are {\em not} just multiples of one another.

In Section \ref{NAS_exist} we give a simple but rigorous proof of the
existence of these non-axisymmetric vortex solutions for all values of
$\lambda \in [0,1)$, provided the circulation Reynolds number of the flow 
is sufficiently small (depending on $\lambda$). A complementary result
is obtained in \cite{gallay:2005} where we prove that, if $\lambda > 0$
is sufficiently small, non-axisymmetric vortex solutions exist for 
{\em all} values of the Reynolds number. 

The construction of these non-axisymmetric vortices requires some
work, so as a first step we rewrite the Navier-Stokes and associated
vorticity equation in non-dimensional form.  This simplifies the
expressions for the solutions and also the equations themselves. Fix
$\lambda \in [0,1)$ and assume that the straining flow is given by
\reff{strain}, \reff{nonsymgam} for some $\gamma > 0$. We replace the
variables $x,t$ and the functions $\uu,p$ with the dimensionless
quantities
\begin{equation}\label{chvar}
  \tilde x \,=\, \Bigl(\frac{\gamma}{\nu}\Bigr)^{1/2}x~, \quad
  \tilde t \,=\, \gamma t~, \quad \tilde \uu \,=\, 
  \frac{\uu}{(\gamma\nu)^{1/2}}~, \quad 
  \tilde p \,=\, \frac{p}{\bar\rho \gamma\nu}~,
\end{equation}
where $\nu > 0$ is the kinematic viscosity. Dropping the tildes
for simplicity, we see that the new functions $\uu,p$ 
satisfy the Navier-Stokes equation \reff{3DNS} with 
$\nu = \bar \rho = 1$. Similarly the new straining flow
$\uu^s$ is given by \reff{strain}, \reff{nonsymgam} with $\gamma = 1$. 
Setting $\uu = \uu^s + \UU$, we obtain that $\OO = \nabla \times 
\UU$ satisfies \reff{3DV1} with $\nu = 1$, namely
\begin{equation}\label{3DV2}
  \partial_t \OO + (\UU\cdot \nabla) \OO - (\OO\cdot\nabla)\UU
   + (\uu^s \cdot \nabla) \OO - (\OO\cdot\nabla)\uu^s
  \,=\, \Delta\OO~, \quad \nabla \cdot \OO \,=\, 0~.
\end{equation}
Thus the main effect of the change of variables \reff{chvar} is to
set $\gamma = \nu = 1$ everywhere. In particular, in the dimensionless 
variables the weight function \reff{rhodef} becomes $b(\xp) = 
(1+x_1^2)^{1/2}(1+x_2^2)^{1/2}$, and the norm \reff{L2mdef} reduces to 
$\|\omega\|_{L^2(m)} = \|b^m \omega\|_{L^2}$. 

To formulate our result, we define
\begin{equation}\label{cGdef}
  {\cG}_{\lambda}(\xp) \,=\, \frac{\sqrt{1-\lambda^2}}{4 \pi} 
  \,e^{-\frac{1}{4}( (1+\lambda)x_1^2+
  (1-\lambda) x_2^2)}~, \quad \xp = (x_1,x_2) \in \real^2~.
\end{equation}

If $\lambda=0$, ${\cG}_\lambda(\xp)$ is just the vorticity
field \reff{hatOBdef} of the symmetric Burgers vortex written in 
the new coordinates, and the family of these vortices is indexed
by the non-dimensionalized circulation number $\rho = \Gamma/\nu$. 
As we show below, for any $\lambda \in (0,1)$, $\cG_{\lambda}(\xp)$ is 
still the leading order approximation to the vorticity of the 
non-axisymmetric Burgers vortex, for small Reynolds number $|\rho|$. 
Our precise result is:

\begin{theorem}\label{asym_exist} Fix $m > 3/2$, $\lambda \in [0,1)$, 
and assume that $(\gamma_1,\gamma_2,\gamma_3)$ is given by 
\reff{nonsymgam} with $\gamma = 1$. There exist $R_1(\lambda) > 0$ 
and $K_1(\lambda) > 0$ such that, for $|\rho| \le R_1 $, the 
vorticity equation \reff{3DV2} has a stationary solution 
$\OO^B(\xp;\rho,\lambda)$ which satisfies
\begin{equation}\label{vortex_circ}
  \OO^B(\xp;\rho,\lambda) \,=\, \pmatrix{0 \cr 0 \cr 
  \Omega^B(\xp;\rho,\lambda)}~, \quad \inttwo 
  \Omega^B(\xp;\rho,\lambda) \d\xp \,=\, \rho~,
\end{equation}
and
\begin{equation}\label{vortex_est}
  \|\Omega^B(\cdot;\rho,\lambda) - \rho {\cG}_{\lambda}
  (\cdot)\|_{L^2(m)} \,\le\, K_1 \rho^2 ~.
\end{equation}
Furthermore, $\Omega^B(\cdot;\rho,\lambda)$ is a smooth function
of $\rho$ and $\lambda$, and there is no other stationary solution 
of \reff{3DV1} of the form \reff{vortex_circ} satisfying 
$\|\Omega^B - \rho {\cG}_{\lambda}\|_{L^2(m)} \le 2R_1$. 
\end{theorem}

\begin{remark}\label{cstrem}
The proof shows that $R_1(\lambda) \to 0$ and $K_1(\lambda) \to 
\infty$ as $\lambda \to 1$. On the other hand, $R_1(0) > 0$ and 
$K_1(\lambda) = \cO(\lambda)$ as $\lambda \to 0$. In particular, 
setting $\lambda = 0$ in \reff{vortex_est}, we recover that
$\Omega^B(\cdot;\rho,0) = \rho {\cG}_{0}$. 
\end{remark}

\begin{remark}\label{smoothrem}
Theorem~\ref{NAS_exist} shows that the asymmetric
Burgers vortex $\OO^B(\xp;\rho,\lambda)$ decays rapidly as 
$|\xp| \to \infty$, since the parameter $m > 3/2$ is arbitrary
(note, however, that the constants $R_1, K_1$ depend on 
$m$). In fact, proceeding as in \cite{gallay:2005}, it is 
possible to show that $\OO^B$ has a Gaussian decay as $|\xp| \to 
\infty$. Moreover, $\OO^B$ is also a smooth function of $\xp$, see 
Remark~\ref{Hkrem} below.
\end{remark}

Finally, we prove that these families of non-symmetric vortices 
are stable in the same sense as the symmetric Burgers vortices are.  

\begin{theorem}\label{asym_stab}
Fix $m > 3/2$, $\lambda \in [0,1)$, and assume that
$(\gamma_1,\gamma_2,\gamma_3)$ is given by \reff{nonsymgam} with
$\gamma = 1$. For any $\mu \in (0,\frac12(1{-}\lambda))$, there 
exist $R_2(\lambda) > 0$ and $\epsilon_2(\lambda) > 0$ such that, 
if $|\rho| \le R_2 $ and if $\OO^0(x) = \OO^B(\xp;\rho,\lambda) 
+ \oomega^0(x)$ with $\oomega^0 \in X^2(m)^3$ satisfying
\begin{equation}\label{phidef}
  \|\oomega^0\|_{X^2(m)} + \lambda
  \|\partial_3 \phi^0\|_{L^\infty}^2 \,\le\, \epsilon_2~,
  \quad \hbox{where} \quad \phi^0(x_3) \,=\, \inttwo
  \omega^0_3(\xp,x_3)\d\xp~,
\end{equation}
then the solution $\OO(x,t)$ of \reff{3DV2} with initial data 
$\OO^0$ converges as $t \to +\infty$ to the vortex solution 
$\OO^B(\xp;\rho',\lambda)$, where $\rho'= \rho + 
\delta \rho$ and $\delta \rho = (2\pi)^{-1/2} \int_\real 
e^{-x_3^2/2}\phi^0(x_3)\d x_3$, see \reff{dRdef}. More precisely, 
for any compact interval $I \subset \real$, we have
\begin{equation}\label{conv2}
  \sup_{x_3 \in I} \|\OO(\cdot,x_3,t)
  - \OO^B(\cdot;\rho',\lambda)\|_{L^2(m)} \,=\, \cO(e^{-\mu t})~,
  \quad t \to +\infty~.
\end{equation}
\end{theorem}

The symmetric case $\lambda=0$ is included in Theorem~\ref{asym_stab},
which therefore subsumes Theorem~\ref{sym_stable}. Note however that
the assumptions on the initial data are more restrictive when 
$\lambda > 0$, because we then need a condition on $\partial_3 
\phi^0$. This is due to the fact that non-axisymmetric Burgers vortices
with different circulation numbers are not multiples of one another.

\medskip
The rest of the text is organized as follows. In Section~\ref{NAS_exist},
we prove the existence of non-axisymmetric Burgers vortices for 
small Reynolds numbers. The core of the paper is Section~\ref{NAS_stab},
where we show that these families of vortices are asymptotically stable
with shift. Section~\ref{appendix} is an appendix where we collect
various estimates on the semigroup associated to the linearized 
vorticity equation, together with a few remarks concerning the 
Biot-Savart law. 


\section{Existence of non-axisymmetric Burgers vortices}
\label{NAS_exist}

The properties of non-axisymmetric Burgers vortices seem first to have
been studied by Robinson and Saffman \cite{robinson:1984} who used
perturbative methods to investigate their existence for small values
of the Reynolds number. There were many further investigations in the
intervening years - we mention particularly the perturbative study of
the large Reynolds number limit of these vortices by Moffatt, Kida and
Ohkitani \cite{moffatt:1994}, and the numerical work of Prochazka and
Pullin \cite{prochazka:1998}. However, as far as we know there has
been no rigorous proof of the existence of these types of solutions
and so in this section we present a simple argument which proves the
existence of non-symmetric vortices in the case of small Reynolds 
number.

Fix $\lambda \in [0,1)$ and assume that $\gamma_1,\gamma_2,\gamma_3$
are given by \reff{nonsymgam} with $\gamma = 1$. Motivated by the 
perturbative calculations of \cite{robinson:1984} we look for 
stationary solutions of \reff{3DV2} of the form
$$
  \OO^B(\xp) \,=\, \pmatrix{0 \cr 0 \cr \Omega^B(\xp)}~, \qquad 
  \inttwo \Omega^B(\xp) \d\xp \,=\, \rho~,
$$
for some $\rho \in \real$ (recall that $|\rho|$ is the Reynolds
number). Since $\OO^B$ depends only on the horizontal variable $\xp =
(x_1,x_2)$ and has only the third component nonzero, the associated
velocity field $\UU^B$ depends only on $\xp$ and has only the first
two components nonzero. Thus $\UU^B$ is naturally identified with a
two-dimensional velocity field $\bar \UU^B$ which can be computed
using the two-dimensional version of the Biot-Savart law:
\begin{equation}\label{BS2D}
  \bar \UU^B(\xp) \,=\, \frac{1}{2\pi} \inttwo
  \frac{1}{|\xp-\yp|^2}\pmatrix{y_2 - x_2 \cr x_1 - y_1}
  \Omega^B(\yp)\d\yp~.
\end{equation}
Inserting these expressions into \reff{3DV2}, we see that 
$\Omega^B$ satisfies the scalar equation
\begin{equation}\label{FP1}
  \bar\UU^B  \cdot \nabla_{\perp} \Omega^B \,=\, 
  (\cLp + \lambda \cM)\Omega^B~,
\end{equation}
where $\cLp$ and $\cM$ are the differential operators
\begin{equation}\label{LMdef}
  \cLp \,=\, \Delta_{\perp} + \half (\xp \cdot \nabla_{\perp}) +1~,
  \quad  \cM \,=\, \frac12 (x_1 \partial_1 - x_2 \partial_2)~.
\end{equation}
Here we have used the natural notations $\nabla_\perp = 
(\partial_1,\partial_2)$ and $\Delta_\perp = 
\partial_1^2 + \partial_2^2$. 

We shall solve \reff{FP1} in the weighted space $L^2(m)$ 
defined by \reff{L2mdef} (with $\gamma = \nu = 1$). 
Our approach rests on the fact that the spectrum of the
linear operator $\cLp + \lambda \cM$ in $L^2(m)$ can be
explicitly computed, see Section~\ref{twodim}. If $m > 1/2$, 
this operator turns out to be invertible on the invariant subspace 
$L^2_0(m)$ defined by 
\begin{equation}\label{L20mdef}
  L^2_0(m) \,=\, \Bigl\{\omega \in L^2(m) \,\Big|\, \inttwo
  \omega(\xp)\d\xp = 0\Bigr\}~.
\end{equation}
This allows to rewrite \reff{FP1} as a fixed point problem which is 
easily solved by a contraction argument. 

As a preliminary step, let $\bar\VV_{\lambda}(\xp)$ be the
two-dimensional velocity field obtained from $\cG_\lambda(\xp)$ by 
the Biot-Savart law \reff{BS2D}. Using \reff{cGdef} and \reff{LMdef} 
one can easily verify that
$$
  (\cLp + \lambda \cM) \cG_{\lambda} = 0~, \quad \hbox{and}~ 
  \inttwo \cG_{\lambda}(\xp) \d\xp = 1~.
$$
If we are given $\Omega^B \in L^2(m)$ with $m > 1/2$ and if
$\rho = \inttwo\Omega^B\d\xp$, we can decompose
\begin{equation}\label{UUOmdecomp}
  \OB \,=\, \rho \cG_{\lambda} + \omega~, \quad 
  \bar\UU^B \,=\, \rho \bar\VV_\lambda + \bar \uu~,
\end{equation}
where $\omega \in L^2_0(m)$ and $\bar \uu$ is the velocity 
obtained from $\omega$ by the Biot-Savart law \reff{BS2D}. 
With these notations, finding a solution to \reff{FP1} is equivalent 
to solving
\begin{equation} \label{FP2}
  (\cLp + \lambda \cM) \omega \,=\, (\rho\bar\VV_{\lambda}+ \bar \uu)
  \cdot \nabla_\perp (\rho \cG_{\lambda} + \omega)~, \quad 
  \omega \in L^2_0(m)~.
\end{equation}
Note that $(\rho\bar\VV_{\lambda}+ \bar \uu)\cdot \nabla_\perp 
(\rho \cG_{\lambda} + \omega) = \nabla_\perp \cdot((\rho
\bar\VV_{\lambda}+ \bar \uu)(\rho\cG_{\lambda} + \omega))$
since $\bar\VV_{\lambda}$ and $\bar\uu$ are divergence-free. 
Thus the right-hand side of the \reff{FP2} has zero mean 
as expected. 

The next proposition ensures that the operator $\cLp + \lambda \cM$ 
is invertible on $L^2_0(m)$ and that $(\cLp + \lambda \cM)^{-1}
\nabla_\perp$ defines a bounded operator from $L^p(m)$ into
$L^2_0(m)$, where $L^p(m)$ is the weighted space defined in 
\reff{Lpmdef}.

\begin{proposition}\label{LMbdd}  
Fix $m > 3/2$ and $\lambda \in [0,1)$. There exists $C(m,\lambda) 
> 0$ such that, for all $f \in L_0^2(m)$, 
\begin{equation}\label{semi1}
  \|(\cLp + \lambda \cM)^{-1}f\|_{L^2(m)} \le C \|f\|_{L^2(m)}~.
\end{equation}
Moreover, if $p \in (1,2]$, there exists $C(m,\lambda,p) > 0$ such 
that, for all $g \in L^p(m)$, 
\begin{equation}\label{semi2}
  \|(\cLp + \lambda \cM)^{-1} \partial_i g\|_{L^2(m)}
  \le C \|g\|_{L^p(m)}~, \quad i=1,2\ .
\end{equation}
\end{proposition}

\proof Let $\cT_{\lambda}(t)$ denote the strongly continuous
semigroup generated by $\cLp + \lambda \cM$, the properties of 
which are studied in the Section~\ref{twodim}. If $f \in L^2_0(m)$, 
we know from \reff{2Dest0} that
\begin{equation}\label{estcT}
  \|\cT_{\lambda}(t)f\|_{L^2(m)} \,\le\, C \,e^{-\frac12(1-\lambda)t}
  \|f\|_{L^2(m)}~, \quad t \ge 0~.
\end{equation}
Thus $\cLp + \lambda \cM$ is invertible on $L^2_0(m)$ and 
we have the Laplace formula
\begin{equation}\label{laplace}
  (\cLp + \lambda \cM)^{-1} f \,=\, -\int_0^{\infty} 
  \cT_{\lambda}(t) f \d t~, \quad f \in L^2_0(m)~.
\end{equation}
Combining \reff{estcT}, \reff{laplace}, we easily obtain
\reff{semi1}. Assume now that $f = \partial_i g$ for some
$i \in \{1,2\}$ and some $g \in L^p(m)$. Using \reff{laplace} 
and \reff{2Dest1} or \reff{2Dest2}, we obtain an estimate of 
the form
$$
  \|(\cLp + \lambda \cM)^{-1} \partial_i g\|_{L^2(m)}
  \,\le\, C \int_0^{\infty} {a(t)^{-\frac{1}{p}}} 
  \,e^{-\half (1-\lambda) t} \|g\|_{L^p(m)}\d t~,
$$
where $a(t) = 1-e^{-t}$. Since $p > 1$, the singularity in the 
integral at $t=0$ is integrable and \reff{semi2} follows.
\QED

\medskip
We can now rewrite \reff{FP2} as $\omega = F_{\lambda,\rho}(\omega)$, 
where $F_{\lambda,\rho} : L^2_0(m) \to L^2_0(m)$ is defined by
\begin{equation}\label{FP3}
  F_{\lambda,\rho}(\omega) \,=\, (\cLp + \lambda \cM)^{-1} 
  \nabla_\perp \cdot \bigl( (\rho\bar\VV_{\lambda} + \bar\uu) \cdot
  (\rho \cG_{\lambda} + \omega) \bigr)~.
\end{equation}
For any $r > 0$, let $B_m(0,r)$ denote the closed ball of radius 
$r$ centered at the origin in $L^2_0(m)$. The main result of this 
section is:

\begin{proposition}\label{fixedpoint}
Fix $m > 3/2$ and $\lambda \in [0,1)$. There exist $R_1(\lambda) 
> 0$ and $K_1(\lambda) > 0$ such that, if $|\rho| \le R_1$, 
then $F_{\lambda,\rho}$ has a unique fixed point $\omega_{\lambda,\rho}$
in $B_m(0,2R_1)$. Moreover $\omega_{\lambda,\rho}$ is contained in 
$B_m(0,K_1\rho^2)$ and $\omega_{\lambda,\rho}$ is a smooth 
function of both $\lambda$ and $\rho$. 
\end{proposition}

\proof Let $\cB_\lambda : L^2(m) \times L^2(m) \to L^2_0(m)$ be the 
bilinear map defined by
$$
  \cB_\lambda(\Omega_1,\Omega_2) \,=\, (\cLp + \lambda \cM)^{-1} 
  \nabla_\perp \cdot (\bar \UU_1 \Omega_2)~,
$$
where $\bar \UU_1$ is the velocity field obtained from $\Omega_1$ by 
the Biot-Savart law \reff{BS2D}. If $1 < p < 2$, then $\|\bar \UU_1 
\Omega_2\|_{L^p(m)} \le C \|\Omega_1\|_{L^2(m)} \|\Omega_2\|_{L^2(m)}$
by Corollary~\ref{2Dprod}. Using in addition \reff{semi2}, we see that
there exists $C_1(\lambda) > 0$ such that
$$
  \|\cB_\lambda(\Omega_1,\Omega_2)\|_{L^2(m)} \,\le\, C_1(\lambda) 
  \|\Omega_1\|_{L^2(m)} \|\Omega_2\|_{L^2(m)}~, \quad 
  \Omega_1, \Omega_2 \in L^2(m)~.
$$
Since $F_{\lambda,\rho}(\omega) = \cB_\lambda(\rho\cG_\lambda +
\omega,\rho\cG_\lambda + \omega)$, we obtain, for all 
$\omega \in L^2_0(m)$, 
\begin{equation}\label{Fmap}
  \|F_{\lambda,\rho}(\omega)\|_{L^2(m)} \,\le\, C_2(\lambda)\rho^2
  + C_3(\lambda)\bigl(2|\rho| \|\omega\|_{L^2(m)} + \|\omega\|_{L^2(m)}^2
  \bigr)~,
\end{equation}
where $C_2(\lambda) = \|\cB_\lambda(\cG_\lambda,\cG_\lambda)\|_{L^2(m)}$
and $C_3(\lambda) = C_1(\lambda)\max(1,\|\cG_\lambda\|_{L^2(m)})$.
Similarly, for all $\omega_1,\omega_2 \in L^2_0(m)$, 
\begin{equation}\label{Fcontract}
  \|F_{\lambda,\rho}(\omega_1) - F_{\lambda,\rho}(\omega_2)\|_{L^2(m)} 
  \,\le\, C_3(\lambda)\|\omega_1 - \omega_2\|_{L^2(m)}
  \bigl(2|\rho| + \|\omega_1\|_{L^2(m)} + \|\omega_2\|_{L^2(m)}\bigr)~.
\end{equation}
When $\lambda = 0$, $\cG_\lambda$ is radially symmetric and
$\bar\VV_\lambda$ is azimuthal, hence $\bar\VV_\lambda \cdot 
\nabla_\perp \cG_\lambda = 0$. Thus $C_2(0) = 0$, so that 
$C_2(\lambda) = \cO(\lambda)$ as $\lambda \to 0$. 

Now, choose $R_1 > 0$ sufficiently small so that
$$
   C_2 R_1 \le 1~, \quad \hbox{and} \quad 8 C_3 R_1 \,\le\, 1~.
$$
If $|\rho| \le R_1$ and $2 C_2 \rho^2 \le r \le 2R_1$, estimates 
\reff{Fmap} and \reff{Fcontract} imply that $F_{\lambda,\rho}$
maps the ball $B_m(0,r)$ into itself and is a strict contraction there. 
More precisely, if $\omega_1, \omega_2 \in B_m(0,r)$, then
$$
  \|F_{\lambda,\rho}(\omega_1)\|_{L^2(m)} \,\le\, r~, \quad \hbox{and} 
  \quad \|F_{\lambda,\rho}(\omega_1) - F_{\lambda,\rho}(\omega_2)\|_{L^2(m)} 
  \,\le\, \frac34 \|\omega_1 - \omega_2\|_{L^2(m)}~. 
$$
By the contraction mapping theorem, $F_{\lambda,\rho}$ has a unique
fixed point $\omega_{\lambda,\rho}$ in $B_m(0,r)$. Choosing $r =
2R_1$, we obtain the existence and uniqueness claim in
Proposition~\ref{fixedpoint}. Then setting $r = K_1 \rho^2$ with $K_1
= 2C_2$, we see that $\omega_{\lambda,\rho} \in B_m(0,K_1 \rho^2)$.
Finally, the smoothness property is a immediate consequence of the
implicit function theorem. Indeed, the map $(\omega,\lambda,\rho)
\mapsto F_{\lambda,\rho}(\omega)$ is obviously $C^\infty$ from
$L^2_0(m) \times [0,1) \times \real$ into $L^2_0(m)$, and the partial
differential
$$
  D_\omega F_{\lambda,\rho}(\omega) \,=\, \tilde \omega \mapsto
  \cB_\lambda(\tilde \omega, \rho\cG_\lambda+\omega) + 
  \cB_\lambda(\rho\cG_\lambda+\omega, \tilde \omega)
$$
satisfies $\|D_\omega F_{\lambda,\rho}(\omega)\| \le 3/4$ 
whenever $|\rho| \le R_1$ and $\omega \in B_m(0,2R_1)$. Thus
$\oone - D_\omega F_{\lambda,\rho}(\omega)$ is invertible at
$\omega = \omega_{\lambda,\rho}$, and the implicit function theorem
implies that $\omega_{\lambda,\rho}$ is a smooth function of both
$\lambda$ and $\rho$. \QED

\medskip
Theorem~\ref{asym_exist} is a direct consequence of 
Proposition~\ref{fixedpoint}. Indeed, if $|\rho| \le R_1(\lambda)$,
we set $\Omega^B(\xp;\rho,\lambda) = \rho \cG_\lambda(\xp) + 
\omega_{\lambda,\rho}(\xp)$, where $\omega_{\lambda,\rho}$ is 
as in Proposition~\ref{fixedpoint}, and we denote by $\bar 
\UU^B(\xp;\rho,\lambda)$ the two-dimensional velocity field obtained
from $\Omega^B$ by the Biot-Savart law \reff{BS2D}. Then
\begin{equation}\label{OmBUB}
  \OO^B(\xp;\rho,\lambda) \,=\, \pmatrix{0 \cr 0 \cr \Omega^B(\xp;
  \rho,\lambda)}~, \quad 
  \UU^B(\xp;\rho,\lambda) \,=\, \pmatrix{\bar U_1^B(\xp;\rho,\lambda) \cr 
  \bar U_2^B(\xp;\rho,\lambda) \cr 0}
\end{equation}
is a stationary solution of \reff{3DV2} which has all the desired
properties. In  particular, since $\omega_{\lambda,\rho} \in L^2_0(m)$,
we have
\begin{equation}\label{integident}
  \inttwo \Omega^B (\xp;\rho,\lambda)\d\xp \,=\, \rho~, 
\end{equation}
while the fact that $\omega_{\lambda,\rho} \in B_m(0,K_1 \rho^2)$
implies that \reff{vortex_est} holds. For later use, we observe 
that there exists $C(\lambda,m) > 0$ such that, for $|\rho| \le R_1$,  
\begin{equation}\label{ombdd}
  \|\Omega^B(\cdot;\rho,\lambda)\|_{L^2(m)} \le C |\rho|~, \quad
  \hbox{and} \quad \|\partial_\rho \Omega^B(\cdot;\rho,\lambda)
  \|_{L^2(m)} \le C~.
\end{equation}
Moreover $\|\partial_\rho^2 \Omega^B(\cdot;\rho,\lambda)
\|_{L^2(m)} \le C\lambda$, because in the symmetric case 
$\Omega^B(\cdot;\rho,0) = \rho \cG_0$ so that $\partial_\rho^2
\Omega^B(\cdot;\rho,0) = 0$. 

\begin{remark}\label{Hkrem}
We chose to solve \reff{FP1} in $L^2(m)$ because this is basically
the space we shall use in Section~\ref{NAS_stab} to study the stability
of the vortices. But it is clear from the proof of 
Proposition~\ref{fixedpoint} that nothing important changes if we 
replace $L^2(m)$ with the corresponding Sobolev space
$$
   H^k(m) \,=\, \Bigl\{f \in L^2(m) \,\Big|\, \partial_1^i 
   \partial_2^j f \in L^2(m) \hbox{ for all }i,j \in \intplus 
   \hbox{ with }i+j \le k\Bigr\}~,
$$
for any $k \in \intplus$. This shows that the asymmetric Burgers
vortex $\OB(\xp;\rho,\lambda)$ is a smooth function of 
$\xp$ too. In particular, by Sobolev embedding, $b^m 
\OB(\cdot;\rho,\lambda) \in C^0_b(\real^2)$ (the space of all 
continuous and bounded functions on $\real^2$) and we have the 
analogue of \reff{ombdd}: 
\begin{equation}\label{ombdd2}
  \sup_{\xp\in\real^2} b(\xp)^m |\Omega^B(\xp;\rho,\lambda)| 
  \le C |\rho|~, \quad
  \sup_{\xp\in\real^2} b(\xp)^m |\partial_\rho \Omega^B(\xp;\rho,\lambda)| 
  \le C~.
\end{equation}
Moreover, since $\Omega^B(\cdot;\rho,\lambda) \in L^p(\real^2)$ for 
all $p \in [1,+\infty]$, Proposition~\ref{BS2Destimates} implies
that $\UU(\cdot;\rho,\lambda) \in L^q(\real^2) \cap C^0_b(\real^2)$ 
for all $q \in (2,\infty]$, and there exists $C(q,m,\lambda) > 0$ 
such that
\begin{equation}\label{ubdd}
  \|\UU^B(\cdot;\rho,\lambda)\|_{L^q(\real^2)} \le C |\rho|~, \quad
  \hbox{and} \quad \|\partial_\rho \UU^B(\cdot;\rho,\lambda)
  \|_{L^q(\real^2)} \le C~.
\end{equation}
\end{remark}


 \section{Stability with respect to three-dimensional 
perturbations}\label{NAS_stab}
 
We now prove that the family of vortices constructed in the previous
section is asymptotically stable with shift, provided the circulation
Reynolds number is sufficiently small, depending on the asymmetry
parameter $\lambda \in [0,1)$. In particular, our result applies to 
the classical family of symmetric Burgers vortices ($\lambda = 0$).

Throughout this section we fix some $\lambda \in [0,1)$. 
For $|\rho|$ sufficiently small we denote by $\OO^B(\xp;\rho)$,
$\UU^B(\xp;\rho)$ the asymmetric vortex \reff{OmBUB} with total
circulation $\rho$ (to simplify the notation, we omit the dependence
on $\lambda$). As we mentioned in the introduction, if we slightly
perturb the vortex $\OO^B(\cdot;\rho)$ the solution of the vorticity
equation will converge toward another vortex with a possibly different
circulation.  This means that we must allow the parameter $\rho$ to
depend on time.  Also, since the perturbations we consider may depend
on the axial variable $x_3$, it turns out to be convenient to
approximate the solutions by vortices with different circulation
numbers in different $x_3$ sections. In other words, we will consider
solutions of \reff{3DV2} of the form
\begin{equation}\label{vorticity_decomp}
  \OO(x,t) \,=\, \pmatrix{0 \cr 0 \cr \OB(\xp;\rho+\phi(x_3,t))}
   \,+\, \pmatrix{\omega_1(x,t) \cr \omega_2(x,t) \cr \omega_3(x,t)}~,
\end{equation}
where $\phi(x_3,t)$ is determined so that $\inttwo\omega_3
(\xp,x_3,t)\d\xp = 0$ for all $x_3$ and $t$. In view of \reff{integident},
it is obvious that any pertubation of $\OO^B(\cdot;\rho)$ that 
is integrable with respect to the transverse variables $\xp$ 
can be decomposed in a unique way as in \reff{vorticity_decomp}.
Similarly, we write the rotational part of the velocity field as 
\begin{equation}\label{velocity_decomp}
  \UU(x,t) = \pmatrix{\tilde U^B_1(x,t;\rho,\phi) \cr 
  \tilde U^B_2(x,t;\rho,\phi) \cr 0} \,+\, \pmatrix{u_1(x,t) \cr 
   u_2(x,t) \cr u_3(x,t)}~,
\end{equation}
where $\tilde \UU^B(x,t;\rho,\phi)$ is the velocity field obtained 
from the vorticity $\OO^B(\xp;\rho+\phi(x_3,t))$ by the Biot-Savart
law \reff{BS3}. It will be shown in Proposition~\ref{UUcompare} that
$\tilde \UU^B(x,t;\rho,\phi)$ is a small perturbation of 
$\UU^B(\xp;\rho+\phi(x_3,t))$ if $\phi$ varies slowly in the 
$x_3$ direction, namely there exists $C(\lambda) > 0$ such 
that
\begin{equation}\label{UtildeU}
  \sup_{x\in\real^3} |\tilde \UU^B(x,t;\rho,\phi) - 
  \UU^B(\xp;\rho+\phi(x_3,t))| \,\le\, C \|\partial_3 
  \phi(\cdot,t)\|_{L^\infty}~.
\end{equation}

Let $\oomega = (\omega_1,\omega_2,\omega_3)^T$ and $\uu = 
(u_1,u_2,u_3)^T$ denote the remainder terms in \reff{vorticity_decomp} 
and \reff{velocity_decomp} respectively. By construction, $\uu$ 
is the velocity field obtained from $\oomega$ by the Biot-Savart
law \reff{BS3}. Remark that $\nabla \cdot \uu = 0$, but $\nabla \cdot 
\oomega = -(\partial_\rho \Omega^B)\partial_3 \phi \neq 0$, hence 
$\oomega \neq \nabla \times \uu$. In broadest terms, our strategy is
to show that $\oomega(x,t)$ and $\partial_3 \phi(x,t)$ converge to
zero as time tends to infinity, so that the vorticity $\OO(x,t)$ 
approaches one of the vortices $\OO^B(\cdot;\rho')$ constructed in 
Section~\ref{NAS_exist}. With that in mind, we now write out the 
evolution equations for $\oomega$ and $\phi$. 

Inserting \reff{vorticity_decomp}, \reff{velocity_decomp} into 
\reff{3DV2} and using the identity $(\UU\cdot\nabla)\OO - 
(\OO\cdot\nabla)\UU = \nabla \times (\OO \times \UU)$, we obtain
after straightforward calculations: 
\begin{equation}\label{o_evolve}
  \partial_t \oomega \,=\, \LLL \oomega + \PPP_\phi \oomega + 
  \NNN(\oomega) + \HHH(\phi)~,
\end{equation}
where the various terms in the right-hand side are defined as follows.
\\[1mm]
$\bullet$ The linear operator $\LLL$ is the leading order part of the 
equation, which takes into account the diffusion and the effects of the
background strain:
$$
  \LLL \oomega \,=\, \Delta \oomega + (\oomega\cdot\nabla)\uu^s 
  - (\uu^s\cdot\nabla)\oomega \,=\, \pmatrix{(\cL + \gamma_1)
  \omega_1 \cr (\cL + \gamma_2)\omega_2 \cr (\cL + \gamma_3)\omega_3}~.
$$
Here $\gamma_1,\gamma_2,\gamma_3$ are given by \reff{nonsymgam} with
$\gamma = 1$, and
\begin{equation}\label{cLdef}
  \cL \,=\, \Delta - (\uu^s\cdot\nabla) \,=\, \Delta + \half 
  (\xp \cdot \nabla_{\perp}) +\frac{\lambda}{2} (x_1 \partial_1 - 
  x_2 \partial_2) - x_3 \partial_3~.
\end{equation}
$\bullet$ The term $\PPP_\phi \oomega = \nabla \times (\tilde \UU^B 
\times \oomega + \uu \times \OO^B)$ describes the linear interaction 
between the perturbation and the modulated vortex, namely:
\begin{equation}\label{PPPexp}
  \PPP_\phi \oomega \,=\, \pmatrix{
  \phantom{-}\partial_2 (\tilde U^B_1 \omega_2 - \tilde U^B_2 \omega_1) +
  \partial_3 (\tilde U^B_1 \omega_3 + u_1 \OB) \cr
  \phantom{-}\partial_1 (\tilde U^B_2 \omega_1 - \tilde U^B_1 \omega_2) +
  \partial_3 (\tilde U^B_2 \omega_3 + u_2 \OB) \cr
  -\partial_1(\tilde U^B_1 \omega_3 + u_1 \OB) 
  -\partial_2(\tilde U^B_2 \omega_3 + u_2 \OB)}~.
\end{equation}
Here and in the sequel, to simplify the notation, we write $\OO^B$ 
instead of $\OO^B(\cdot;\rho+\phi)$ and $\tilde \UU^B$ instead 
of $\tilde\UU^B(\cdot;\rho,\phi)$.
\\[1mm]
$\bullet$ The term $\NNN(\oomega) = \nabla \times (\uu \times \oomega)$ 
collects all the nonlinear contributions in $\oomega$, specifically:
\begin{equation}\label{NNNexp}
  \NNN(\oomega) \,=\, \pmatrix{
  \partial_2 (u_1\omega_2 - u_2\omega_1) + 
  \partial_3 (u_1\omega_3 - u_3\omega_1) \cr
  \partial_1 (u_2\omega_1 - u_1\omega_2) + 
  \partial_3 (u_2\omega_3 - u_3\omega_2) \cr
  \partial_1 (u_3\omega_1 - u_1\omega_3) + 
  \partial_2 (u_3\omega_2 - u_2\omega_3)}~.
\end{equation}
$\bullet$ Finally, $\HHH(\phi) = \LLL \OO^B + \nabla \times(\tilde\UU^B 
\times \OO^B) - \partial_t \OO^B$ is an inhomogeneous term which is 
due to the fact that $\OO^B(\cdot;\rho+\phi)$ fails to be a solution of 
\reff{3DV2} if $\phi$ is not identically constant. A simple
calculation shows that $\HHH_i(\phi) = \partial_3(\tilde U_i^B 
\Omega^B)$ for $i = 1,2$. The third component of $\HHH(\phi)$ has
a more complicated expression:
\begin{eqnarray}\nonumber
  \HHH_3(\phi) &=& (\cLp + \lambda\cM)\Omega^B - \nabla_\perp
  \cdot(\tilde \UU^B \Omega^B) + (\partial_{\rho}^2 \OB)(\partial_3
  \phi)^2 \\ \label{HHH3def}
  && -(\partial_\rho \OB)(\partial_t \phi + x_3 \partial_3 \phi 
  -\partial_3^2 \phi)~,
\end{eqnarray}
where $\cLp$ and $\cM$ are defined in \reff{LMdef}. 

Equation \reff{o_evolve} governs the evolution of both $\phi$ and
$\oomega$. To separate out the evolution equation for $\phi$, we
recall that $\oomega$ satisfies the constraint $\inttwo \omega_3
(\xp,x_3,t)\d\xp = 0$. If we integrate the third component of the
vectorial equation \reff{o_evolve} with respect to the transverse
variables $\xp$, the first three terms in the right-hand side give no
contribution, as can be seen from the formulas \reff{cLdef},
\reff{PPPexp}, \reff{NNNexp}. So we must impose
\begin{equation}\label{HHHconst}
  \inttwo \HHH_3(\phi)\d\xp \,=\, 0~, \quad \hbox{for all } 
  x_3 \hbox{ and } t~.
\end{equation}
As is clear from \reff{FP1}, the first term in the right-hand side
of \reff{HHH3def} has zero mean with respect to $\xp$, and so does
the second term because it is explicitly in divergence form. 
On the other hand, differentiating \reff{integident} with respect 
to $\rho$, we obtain the identities $\inttwo \partial_\rho \OB\d\xp 
= 1$ and $\inttwo\partial_\rho^2 \OB\d\xp = 0$. Thus \reff{HHHconst} 
gives the evolution equation for $\phi$:
\begin{equation}\label{phi_evolve}
  \partial_t \phi + x_3 \partial_3 \phi \,=\, \partial_3^2 \phi~.
\end{equation}

Remarkably, this equation is {\em linear} and completely {\em 
decoupled} from the rest of the system. As is easily verified, the 
solution of \reff{phi_evolve} 
with initial data $\phi(x_3,0)=\phi^0(x_3)$ is given by the explicit
formula
\begin{equation}\label{phi_solve}
  \phi(x_3,t) \,=\, (G_t * \phi^0)(x_3 e^{-t})~, \quad x_3 \in \real~, 
  \quad t > 0~,
\end{equation}
where 
\begin{equation}\label{Gtdef}
  G_t(x_3) \,=\, \sqrt{\frac{1}{2\pi(1{-}e^{-2 t})}}
  \,\exp\Bigl({-\frac{x_3^2}{2(1{-}e^{-2 t})}}\Bigr)~, \quad x_3 
  \in \real~, \quad t > 0~.  
\end{equation}
The following simple estimates will be useful:

\begin{proposition}\label{phi_estimates} If $\phi^0 \in C^0_b(\real)$, 
the solution of \reff{phi_evolve} with initial data $\phi^0$ 
satisfies
\begin{equation}\label{phiest1}
  \|\phi(\cdot,t)\|_{L^\infty} \,\le\, \|\phi^0\|_{L^\infty}~, \quad
  \|\partial_3 \phi(\cdot,t)\|_{L^\infty} \,\le\, \frac{e^{-t}}
  {\sqrt{1-e^{-2t}}}\|\phi^0\|_{L^\infty}~, \quad t > 0~.
\end{equation}
If moreover $\partial_3 \phi^0 \in L^\infty(\real)$, we also have
\begin{equation}\label{phiest2}
  \|\partial_3 \phi(\cdot,t)\|_{L^\infty} \,\le\, e^{-t}
  \|\partial_3 \phi^0\|_{L^\infty}~, \quad t \ge 0~.
\end{equation}
\end{proposition}

\proof
Since $\|G_t\|_{L^1} = 1$, it follows immediately from \reff{phi_solve}
that $\|\phi(\cdot,t)\|_{L^\infty} \le\|\phi^0\|_{L^\infty}$. If 
$\partial_3 \phi^0 \in L^\infty(\real)$, the same argument gives 
\reff{phiest2}, because
\begin{equation}\label{dphiid}
  \partial_3 \phi(x_3,t) \,=\, e^{-t}(G_t*\partial_3\phi^0)(x_3 e^{-t})
  \,=\, e^{-t}(\partial_3 G_t * \phi^0)(x_3 e^{-t})~, \quad t > 0~.
\end{equation}
To prove the second estimate in \reff{phiest1}, we use the last 
expression in \reff{dphiid} and observe that $\|\partial_3 G_t
\|_{L^1} = C/\sqrt{1-e^{-2t}}$, where $C = \sqrt{2/\pi} < 1$. \QED

\begin{remark}\label{notcontinuous}
Proposition~\ref{phi_estimates} shows in particular that 
\reff{phi_solve} defines a semigroup of bounded linear operators
on $C^0_b(\real)$, the space of all bounded and continuous functions
on $\real$ equipped with the $L^\infty$ norm. It is easy to verify 
that this semigroup is {\em not strongly continuous} in time, due
to the dilation factor $e^{-t}$ in \reff{phi_solve} which in turn 
originates in the unbounded advection term $x_3\partial_3$ in 
\reff{phi_evolve}. However, if we equip $C^0_b(\real)$ with the 
(weaker) topology of uniform convergence on compact sets, then
\reff{phi_solve} defines a continuous function of time. 
\end{remark}

We now return to the evolution equation for $\oomega$. Using 
\reff{HHH3def}, equation \reff{phi_evolve} for $\phi$, and 
equation \reff{FP1} satisfied by the asymmetric vortex $\Omega^B$, 
we obtain for the inhomogeneous term $\HHH(\phi)$ the simpler 
expression
\begin{equation}\label{HHHexp}
  \HHH(\phi) \,=\, \pmatrix{
  \partial_3 (\tilde U^B_1 \OB)\cr
  \partial_3 (\tilde U^B_2 \OB)\cr
  \nabla_\perp \cdot ((\UU^B{-}\tilde \UU^B)\Omega^B) + 
  (\partial_{\rho}^2 \OB) (\partial_3 \phi)^2}~,
\end{equation}
where as usual $\UU^B = \UU^B(\cdot;\rho+\phi)$. Before starting the
rigorous analysis, let us briefly comment here on why we expect
solutions of \reff{o_evolve} to go to zero as $t$ goes to infinity.
Given $m > 3/2$, we assume that $\omega_i \in X^2(m)$ for $i = 1,2,3$,
where $X^2(m)$ is the space defined in \reff{X2mdef}. By construction,
$\omega_3$ then belongs to the subspace $X^2_0(m)$ given by
\reff{X20mdef}. As we show in Section~\ref{threedim}, the linear
operator $\LLL$ has spectrum that lies in the half-plane $\{z \in
\complex ~|~ \Re z \le -\half(1{-}\lambda)\}$ when acting on $X^2(m)
\times X^2(m) \times X^2_0(m)$. Thus, the semigroup generated by this
operator can be expected to decay like $\exp(-\half(1{-} \lambda)t)$.
The remaining linear terms in the equation, namely
$\PPP_\phi(\omega)$, contain a factor of the vortex solution which is
proportional to $\rho+\phi$ (see \reff{ombdd2} and \reff{UUfirst}) and
hence, for small Reynolds number, they will be a small perturbation of
$\LLL$ and will not destroy the exponential decay. The same is true
for the nonlinear terms $\NNN(\omega)$ provided we restrict
ourselves to sufficiently small perturbations.  Finally, the
inhomogeneous term $\HHH(\phi)$ decays at least like $e^{-t}$ by
\reff{UtildeU} and Proposition~\ref{phi_estimates}, so we expect the
solution $\oomega(x,t)$ of \reff{o_evolve} to converge exponentially
to zero if the initial data are sufficiently small. 

We now put these heuristic arguments into a rigorous form.  Let
$\XXX(m)$ be the Banach space $X^2(m) \times X^2(m) \times X^2_0(m)$
equipped with the norm $\|\oomega\|_{\XXX(m)} = \|\omega_1\|_{X^2(m)}
+ \|\omega_2\|_{X^2(m)} + \|\omega_3\|_{X^2(m)}$. As is shown in
Proposition~\ref{3Dsemig}, the linear operator $\cL$ is the generator 
of a semigroup $e^{t\cL}$ of bounded operators on $X^2(m)$, hence
the same is true for the operator $\LLL$ acting on $\XXX(m)$. A
natural idea is then to use Duhamel's formula to rewrite
\reff{o_evolve} as an integral equation:
\begin{equation}\label{duhamel}
  \oomega(t) \,=\, e^{t \LLL}\oomega^0 + \int_0^t e^{(t-s) \LLL} 
  \Bigl(\PPP_\phi\oomega(s) + \NNN(\oomega(s)) + \HHH(\phi(s))\Bigr) 
  \d s~, \quad t \ge 0~,
\end{equation}
which can then be solved by a fixed point argument. A problem with 
this approach is that the semigroup $e^{t\cL}$ fails to be strongly
continuous on $X^2(m)$, essentially for the reason mentioned in 
Remark~\ref{notcontinuous}. To restore continuity, it is thus necessary
to equip $X^2(m)$ with a weaker topology. For any $n \in \intplus^*$ 
we define the seminorm
\begin{equation}\label{2seminorm}
  |\omega|_{X^2_n(m)} \,=\, \sup_{|x_3|\le n} 
  \|\omega(\cdot,x_3)\|_{L^2(m)}~,
\end{equation}
and we denote by $X^2_\loc(m)$ the space $X^2(m)$ equipped with 
the topology defined by the family of seminorms \reff{2seminorm}
for $n \in \intplus^*$, i.e. the topology of the Fr\'echet space 
$C^0(\real,L^2(m))$. In other words, a sequence $\omega_k$ converges
to zero in $X^2_\loc(m)$ if and only if $|\omega_k|_{X^2_n(m)} \to 0$
as $k \to \infty$ for all $n \in \intplus^*$, namely if
$\omega_k(x_3)$ converges to zero in $L^2(m)$ {\em uniformly on compact 
sets} in $x_3$. We define the product space $\XXX_\loc(m)$ in a 
similar way. Then Proposition~\ref{3Dsemig} shows that the semigroup
$e^{t\LLL}$ is strongly continuous on $\XXX_\loc(m)$, and the integrals
in \reff{duhamel} can be defined as $\XXX_\loc(m)$-valued Riemann 
integrals, see Corollary~\ref{3Dint} and Remark~\ref{3Dint2}. 

Since we expect $\oomega(t)$ to converge exponentially to zero
as $t \to \infty$, we shall solve \reff{duhamel} in the Banach
space
$$
  \YYY_\mu(m) \,=\, \{\oomega \in C^0([0,+\infty),\XXX_\loc(m))
  \,|\, \|\oomega\|_{\YYY_\mu(m)} < \infty\}~,
$$
for some $\mu > 0$, where
$$
  \|\oomega\|_{\YYY_\mu(m)} \,=\, \sup_{t \ge 0} e^{\mu t}
  \|\oomega(t)\|_{\XXX(m)}~.
$$ 
Given initial data $\phi^0 \in C^0_b(\real)$ and $\oomega^0 \in
\XXX(m)$, we first define $\phi(x_3,t)$ by \reff{phi_solve}, and then
use the integral equation \reff{duhamel} to determine $\oomega(t)$ 
for all $t \ge 0$. Our main result is:

\begin{proposition}\label{exproof}
Fix $\lambda \in [0,1)$, $m > 3/2$, and $0 < \mu < \half(1{-}\lambda)$. 
There exist positive constants $\rho_2 > 0$, $\epsilon_2 > 0$, and 
$K_2 > 0$ such that, if $|\rho| \le \rho_2$, $\epsilon \le \epsilon_2$,
and if $\phi^0 \in C^0_b(\real)$ satisfies $\|\phi^0\|_{L^\infty} + 
\lambda \|\partial_3 \phi^0\|_{L^\infty}^2 \le \epsilon$, then for all
$\oomega^0 \in \XXX(m)$ with $\|\oomega^0\|_{\XXX(m)} \le \epsilon$ 
equation \reff{duhamel} has a unique solution $\oomega \in \YYY_\mu(m)$
with $\|\oomega\|_{Y_\mu(m)} \le K_2 \epsilon$. 
\end{proposition}

\proof
Fix $\lambda \in [0,1)$, $m > 3/2$, and $0 < \mu < \half(1{-}\lambda)$. 
To simplify the notations, we shall write $\XXX$ instead of $\XXX(m)$
and $\YYY$ instead of $\YYY_\mu(m)$. For any $r > 0$, we denote by
$B_X(0,r)$ (respectively, $B_Y(0,r)$) the closed ball of radius 
$r > 0$ centered at the origin in $\XXX$ (respectively, $\YYY$).
Let $\phi^0 \in C^0_b(\real)$ and denote by $\phi(x_3,t)$ the 
solution of \reff{phi_evolve} with initial data $\phi^0$. 
Given $\rho \in \real$, $\oomega^0 \in \XXX$, and $\oomega \in \YYY$, 
we estimate the various terms in the right-hand side of \reff{duhamel}. 

We begin with the linear term $e^{t\LLL}\oomega^0$. From 
Proposition~\ref{3Dsemig} we know that the linear operator $\cL + 1
\equiv \hat \cL_{1+\lambda,1-\lambda}$ generates a semigroup 
$S_t = e^{t(\cL+1)}$ which is strongly continuous on $X^2_\loc(m)$, 
uniformly bounded on $X^2(m)$, and which decays like 
$e^{-\half(1{-}\lambda)t}$ on $X^2_0(m)$. Since
$$
  e^{t\LLL}\oomega^0 \,=\, \Bigl(e^{t(\cL+\gamma_1)}\omega^0_1\,,\, 
  e^{t(\cL+\gamma_2)}\omega^0_2\,,\, e^{t(\cL+\gamma_3)}\omega^0_3\Bigr)^T~,
$$
where $\gamma_1,\gamma_2,\gamma_3$ are given by \reff{nonsymgam}
with $\gamma = 1$, we deduce that $t \mapsto e^{t\LLL}\oomega^0$
is continuous in $\XXX_\loc$ and satisfies
\begin{eqnarray}\nonumber
  \|e^{t\LLL}\oomega^0\|_\XXX &\le& C\Bigl(e^{-\frac{3+\lambda}
  {2}t}\|\omega^0_1\|_{X^2(m)} + e^{-\frac{3-\lambda}
  {2}t}\|\omega^0_2\|_{X^2(m)} + e^{-\frac{1-\lambda}{2}t}
  \|\omega^0_3\|_{X^2(m)}\Bigr)\\ \label{LLLest} 
  &\le& C_1 e^{-\frac{1-\lambda}{2}t}\|\oomega^0\|_\XXX~.
\end{eqnarray} 
Note that it is crucial here that $\omega^0_3 \in X^2_0(m)$, otherwise
we do not get any decay at all. 

We next consider the linear term $\int_0^t e^{(t-s)\LLL}\PPP_\phi
\oomega(s)\d s$. For $s \ge 0$ and $i \in \{1,2\}$, we know from 
Proposition~\ref{UUcompare} that $\tilde U_i^B(s) \equiv \tilde
 U_i^B(\cdot,s;\rho,\phi(s)) \in C^0_b(\real^3)$ and 
$\|\tilde U_i^B(s)\|_{L^\infty} \le C(|\rho|+\|\phi(s)\|_{L^\infty})
\le C(|\rho|+\|\phi^0\|_{L^\infty})$ by Proposition~\ref{phi_estimates}. 
Since $\omega_j(s) \in X^2(m)$ for $j \in \{1,2,3\}$, it follows that
$\tilde U_i^B(s)\omega_j(s) \in X^2(m)$ and
$$
  \|\tilde U_i^B(s)\omega_j(s)\|_{X^2(m)} \,\le\, 
  \|\tilde U_i^B(s)\|_{L^\infty(\real^3)} \|\omega_j(s)\|_{X^2(m)} 
  \,\le\, C\rho' \|\omega_j(s)\|_{X^2(m)}~,
$$
where $\rho' = |\rho|+\|\phi^0\|_{L^\infty}$.  Moreover, it is not
difficult to verify that $s \mapsto \tilde U_i^B(s) \omega_j(s)$ is
continuous in $X^2_\loc(m)$. Similarly, $u_i(s) \in X^q(0)$ for $q \in
(2,+\infty)$ by Proposition~\ref{velocity-vorticity} and $\Omega^B(s)
\equiv \Omega^B(\cdot;\rho+\phi(s)) \in X^p(m)$ for $p \in
[1,+\infty]$ by Remark~\ref{Hkrem}, hence $u_i(s)\Omega^B(s) \in
X^2(m)$ by H\"older's inequality and
$$
  \|u_i(s)\Omega^B(s)\|_{X^2(m)} \,\le\, 
  \|u_i(s)\|_{X^4(0)} \|\Omega^B(s)\|_{X^4(m)} 
  \,\le\, C\rho' \|\omega_i(s)\|_{X^2(m)}~.
$$
Again $s \mapsto u_i(s)\Omega^B(s)$ is continuous in $X^2_\loc(m)$.
Thus Corollary~\ref{3Dint} and Remark~\ref{3Dint2} show that the
three components of the vector $\int_0^t e^{(t-s)\LLL}\PPP_\phi
\oomega(s)\d s$ are well defined and continuous in $X^2_\loc(m)$ 
for $t \ge 0$. Using Proposition~\ref{3Dsemig}, we can estimate the 
first component as follows:
\begin{eqnarray*}
  &&\Bigl\|\int_0^t e^{(t-s)(\cL+\gamma_1)} \Bigl(\partial_2(\tilde U^B_1 
  \omega_2 - \tilde U^B_2 \omega_1) + \partial_3 (\tilde U^B_1
  \omega_3 + u_1 \OB)\Bigr)(s)\d s \Bigr\|_{X^2(m)} \\ 
  &&\quad \le C\int_0^t \frac{e^{-2(t-s)}}{a(t{-}s)^{1/2}} 
  \Bigl(\|\tilde U^B_1(s)\omega_2(s)\|_{X^2(m)} + \|\tilde U^B_2(s)
  \omega_1(s)\|_{X^2(m)}\Bigr)\d s \\ 
  &&\quad  +~ C\int_0^t \frac{e^{-\frac{3+\lambda}{2}(t-s)}}{a(t{-}s)^{1/2}} 
  \Bigl(\|\tilde U^B_1(s)\omega_3(s)\|_{X^2(m)} + \|u_1(s)
  \Omega^B(s)\|_{X^2(m)}\Bigr)\d s~,
\end{eqnarray*}
where $a(t) = 1 - e^{-t}$. (This estimate could be sharpened somewhat
by using the functions $a_1(t)$, $a_2(t)$, and $1-e^{-2t}$ which appear
in Proposition \ref{3Dsemig}, but they would lead to no qualitative
improvement in the final result and so we use this somewhat simpler
form.) The other two components can be estimated in exactly the same 
way except for a slower exponential decay of the linear semigroup, 
see \reff{LLLest}. Summarizing, we obtain:
\begin{equation}\label{PPPest}
  \Bigl\|\int_0^t e^{(t-s)\LLL}\PPP_\phi\oomega(s)\d s\Bigr\|_{\XXX}
  \,\le\, C\rho' \int_0^t \frac{e^{-\frac{1-\lambda}{2}(t-s)}}
  {a(t{-}s)^{1/2}} \,\|\oomega(s)\|_\XXX \d s \,\le\, C_2\rho'
  e^{-\mu t} \|\oomega\|_\YYY~.
\end{equation}

We now consider the nonlinear term $\int_0^t e^{(t-s)\LLL}\NNN(\oomega(s))
\d s$. Let $1 < p < 2$. For $s \ge 0$ and $i,j \in \{1,2,3\}$, we know
from Corollary~\ref{3Dprod} that $u_i(s)\omega_j(s) \in X^p(m)$ with
$$
  \|u_i(s)\omega_j(s)\|_{X^p(m)} \,\le\, C \|\omega_i(s)\|_{X^2(m)}
  \|\omega_j(s)\|_{X^2(m)}~.
$$
Moreover $s \mapsto u_i(s)\omega_j(s)$ is continuous in $X^p_\loc(m)$. 
Thus, by Remark~\ref{3Dint2}, the integral $\int_0^t e^{(t-s)\LLL}
\NNN(\oomega(s))\d s$ is well defined and continuous in $\XXX_\loc$
for $t \ge 0$. Proceeding as above we can estimate the first 
component as follows:
\begin{eqnarray*}
  &&\Bigl\|\int_0^t e^{(t-s)(\cL+\gamma_1)} \Bigl(\partial_2 
  (u_1\omega_2 - u_2\omega_1) + \partial_3 (u_1\omega_3 - u_3\omega_1)
  \Bigr)(s)\d s \Bigr\|_{X^2(m)} \\
  &&\quad \le C\int_0^t \frac{e^{-2(t-s)}}{a(t{-}s)^{1/p}} 
  \Bigl(\|u_1(s)\omega_2(s)\|_{X^p(m)} + \|u_2(s)\omega_1(s)\|_{X^p(m)}
  \Bigr)\d s \\
  &&\quad  +~C\int_0^t \frac{e^{-\frac{3+\lambda}{2}(t-s)}}{a(t{-}s)^{1/p}} 
  \Bigl(\|u_1(s)\omega_3(s)\|_{X^p(m)} + \|u_3(s)\omega_1(s)\|_{X^p(m)}
  \Bigr)\d s~.
\end{eqnarray*}
Repeating the same arguments for the other two components, we thus 
find
\begin{equation}\label{NNNest}
  \Bigl\|\int_0^t e^{(t-s)\LLL}\NNN(\oomega(s))\d s\Bigr\|_{\XXX}
  \le C \int_0^t \frac{e^{-\frac{1-\lambda}{2}(t-s)}}{a(t{-}s)^{1/p}} 
  \,\|\oomega(s)\|_\XXX^2 \d s \,\le\, C_3 e^{-\mu t} \|\oomega\|_\YYY^2~.
\end{equation}
 
Finally, we turn our attention to the inhomogeneous term $\int_0^t 
e^{(t-s)\LLL}\HHH(\phi(s))\d s$. To bound the first two components, we 
observe that $\partial_3(\tilde U_i^B(s) \Omega^B(s)) \in X^2(m)$
with
\begin{equation}\label{HHH1}
  \|\partial_3(\tilde U_i^B(s) \Omega^B(s))\|_{X^2(m)} \,\le\, 
  C \rho' \|\partial_3 \phi(s)\|_{L^\infty}~,
  \quad i = 1,2.
\end{equation}
Indeed, $\partial_3(\tilde U_i^B \Omega^B) = \tilde U_i^B
\partial_3\Omega^B + \Omega^B \partial_3 \tilde U_i^B$. Since
$\Omega^B = \Omega^B (\xp;\rho+\phi(x_3,t))$, one has $\partial_3
\Omega^B = (\partial_\rho \Omega^B)(\partial_3 \phi)$, hence by
\reff{ombdd}, \reff{UUfirst}
$$
   \|\tilde U_i^B\partial_3\Omega^B \|_{X^2(m)}\,\le\, 
   \|\tilde U_i^B\|_{L^\infty}\|\partial_3\Omega^B\|_{X^2(m)}
   \,\le\, C \rho' \|\partial_3 \phi\|_{L^\infty}~.
$$
Moreover, as $\partial_3 \tilde \UU^B$ is the velocity field obtained 
from $\partial_3 \OO^B = (\partial_\rho \OO^B)(\partial_3\phi)$ 
by the Biot-Savart law, the proof of Proposition~\ref{UUcompare} 
shows that $\|\partial_3 \tilde U_i^B\|_{L^\infty} 
\le C\|\partial_3\phi \|_{L^\infty}$, hence
$$
   \|\Omega^B \partial_3 \tilde U_i^B\|_{X^2(m)}\,\le\, 
   \|\Omega^B\|_{X^2(m)} \|\partial_3 \tilde U_i^B\|_{L^\infty}
   \,\le\, C \rho' \|\partial_3 \phi\|_{L^\infty}~,
$$
which proves \reff{HHH1}. As usual, one checks that $s \mapsto
\partial_3(\tilde U_i^B(s)\Omega^B(s))$ is continuous in 
$X^2_\loc(m)$ for $s > 0$. Using \reff{phiest1}, \reff{HHH1}, the
first component of the inhomogeneous term can be estimated as 
follows:
\begin{eqnarray*}
  \Bigl\|\int_0^t e^{(t-s)(\cL+\gamma_1)} \partial_3(\tilde U_1^B 
  \Omega^B)(s)\d s \Bigr\|_{X^2(m)} 
  &\le& C\rho' \int_0^t e^{-\frac{3+\lambda}{2}(t-s)}\frac{e^{-s}}
  {a(s)^{1/2}}\|\phi^0\|_{L^\infty} \d s \\
  &\le& C e^{-\mu t}\rho' \|\phi^0\|_{L^\infty}~,
\end{eqnarray*}
and the second one is bounded in exactly the same way. To bound the
third component, we first remark that $(U_i^B(s)-\tilde U_i^B(s))
\Omega^B(s)$ belongs to $X^2(m)$ for $i = 1,2$ and depends 
continuously on $s > 0$ in $X^2_\loc(m)$. By \reff{ombdd}, 
\reff{UtildeU}, \reff{phiest1}, 
$$
  \|(U_i^B-\tilde U_i^B)\Omega^B\|_{X^2(m)} \,\le\,
  \|U_i^B-\tilde U_i^B\|_{L^\infty}\|\Omega^B\|_{X^2(m)}
  \,\le\, C \rho' \|\partial_3 \phi\|_{L^\infty}~,
$$
hence using \reff{phiest1} we find
\begin{eqnarray*}
  &&\Bigl\|\int_0^t e^{(t-s)(\cL+\gamma_3)} \nabla_\perp \cdot 
  (\UU^B(s)-\tilde\UU^B(s))\Omega^B(s)\d s \Bigr\|_{X^2(m)} \\
  && \qquad \qquad \,\le\, C\rho' \int_0^t 
  \frac{e^{-\frac{1-\lambda}{2}(t-s)}}
  {a(t-s)^{1/2}}\frac{e^{-s}}{a(s)^{1/2}}\|\phi^0\|_{L^\infty} \d s
  \,\le\, C e^{-\mu t}\rho' \|\phi^0\|_{L^\infty}~.
\end{eqnarray*}
On the other hand, $\partial_\rho^2 \Omega^B(s) (\partial_3 \phi(s))^2$
lies in $X^2(m)$ and depends continuously on $s > 0$ in $X^2_\loc(m)$.
As was mentioned before Remark~\ref{Hkrem}, $\|\partial_\rho^2 
\Omega^B(s)\|_{X^2(m)} \le C\lambda$. If $\lambda > 0$, we assume
that $\partial_3 \phi^0 \in L^\infty(\real)$ and using \reff{phiest2}
we obtain
\begin{eqnarray*}
  \Bigl\|\int_0^t e^{(t-s)(\cL+\gamma_3)} (\partial_\rho^2 \Omega^B(s)) 
  (\partial_3\phi(s))^2 \d s \Bigr\|_{X^2(m)} &\le& 
  C\lambda \int_0^t e^{-\frac{1-\lambda}{2}(t-s)} 
  e^{-2s}\|\partial_3 \phi^0\|_{L^\infty}^2 \d s \\ &\le&
  C \lambda e^{-\mu t}\|\partial_3 \phi^0\|_{L^\infty}^2~.
\end{eqnarray*}
Thus we have shown that $\int_0^t e^{(t-s)\LLL}\HHH(\phi(s))\d s$ is 
well defined and continuous in $\XXX_\loc$ for $t \ge 0$. 
Moreover, 
\begin{equation}\label{HHHest}
  \Bigl\|\int_0^t e^{(t-s)\LLL}\HHH(\phi(s))\d s\Bigr\|_{\XXX}
  \le C_4 e^{-\mu t}\rho' \|\phi^0\|_{L^\infty} + C_5 \lambda 
  e^{-\mu t}\|\partial_3 \phi^0\|_{L^\infty}^2~.
\end{equation}

Given $\oomega \in \YYY$, we denote by $(\FFF\oomega)(t)$ the right-hand 
side of \reff{duhamel}. Estimates \reff{LLLest}, \reff{PPPest}, 
\reff{NNNest}, \reff{HHHest} show that $t \mapsto (\FFF\oomega)(t)
\in \YYY$. Thus $\FFF$ maps $\YYY$ into itself and
\begin{equation}\label{FFFmap}
  \|\FFF\oomega\|_\YYY \,\le\, C_1 \|\oomega^0\|_\XXX + C_2 \rho'
  \|\oomega\|_\YYY + C_3 \|\oomega\|_\YYY^2 + C_4 \rho' 
  \|\phi^0\|_{L^\infty} + C_5 \lambda \|\partial_3
  \phi^0\|_{L^\infty}^2~,
\end{equation}
where $\rho' = |\rho| + \|\phi^0\|_{L^\infty}$. Moreover, if 
$\oomega_1, \oomega_2 \in \YYY$, the same estimates show that
\begin{equation}\label{FFFcont}
  \|\FFF\oomega_1 -\FFF\oomega_2\|_\YYY \,\le\, 
  \|\oomega_1 - \oomega_2\|_\YYY \Bigl(C_2 \rho'
  + C_3 (\|\oomega_1\|_\YYY + \|\oomega_2\|_\YYY)\Bigr)~,
\end{equation}
because the linear term $e^{t\LLL}\oomega^0$ and the inhomogeneous
term depending on $\HHH(\phi)$ drop out when we consider the difference 
$\FFF\oomega_1 -\FFF\oomega_2$. Now, choose $\rho_2 > 0$ and 
$\epsilon_2 > 0$ small enough so that
$$
  \rho_2 + \epsilon_2 \,\le\, \min\Bigl(1,R_1,\frac{1}{2C_2}\Bigr)~, 
  \quad \hbox{and}\quad \epsilon_2 \le \frac{1}{32C_3
  (C_1+C_4+C_5)}~,
$$
where $R_1$ is as in Proposition~\ref{fixedpoint}. 
Assume that $|\rho| \le \rho_2$, $\epsilon \le \epsilon_2$,
$\|\phi^0\|_{L^\infty} + \lambda \|\partial_3 \phi^0\|_{L^\infty}^2 
\le \epsilon$, and $\|\oomega^0\|_\XXX \le \epsilon$. 
If $4(C_1+C_4+C_5)\epsilon \le r \le 1/(8C_3)$, then \reff{FFFmap} 
shows that $\FFF$ maps the ball $B_Y(0,r)$ into itself. Indeed, 
under the assumptions above we have $C_2 \rho' \le 1/2$ and
$C_3 r \le 1/4$, hence if $\oomega \in B_Y(0,r)$ then \reff{FFFmap} 
implies
$$
  \|\FFF\oomega\|_\YYY \,\le\, C_1 \epsilon + \frac{r}2 + 
  \frac{r}{4} + C_4\epsilon + C_5\epsilon \,=\,
  (C_1+C_4+C_5)\epsilon + \frac{3r}4 \,\le\, r~.
$$
Similarly, $\|\FFF\oomega_1 -\FFF\oomega_2\|_\YYY \le \frac34
\|\oomega_1 - \oomega_2\|_\YYY$ if $\oomega_1, \oomega_2 \in 
B_Y(0,r)$. By the contraction mapping theorem, $\FFF$ has 
thus a unique fixed point $\oomega$ in $B_Y(0,r)$. 
Choosing $r = K_2 \epsilon$ with $K_2 = 4(C_1+C_4+C_5)$, we see
that $\oomega$ is the unique solution of \reff{duhamel} 
such that $\|\oomega\|_\YYY \le K_2\epsilon$.\QED

\medskip Theorem~\ref{asym_stab} is a direct consequence of
Proposition~\ref{exproof}. Indeed, suppose that the initial condition
for the vorticity is $\OO^0(x) = \OO^B(\xp;\rho) + \oomega^0(x)$,
where $\oomega^0 \in X^2(m)^3$ satisfies \reff{phidef}.
Then we can decompose $\OO^0(x) = \OO^B(\xp;\rho+\phi^0(x_3)) +
\tilde\oomega^0(x)$, where $\phi^0$ is as in \reff{phidef} and
$\tilde\oomega^0$ belongs to $\XXX(m)$, namely $\tilde\omega^0_3 \in
X^2_0(m)$. Moreover, there exists $C(m,\lambda) > 0$ such that
$$
  \|\tilde \oomega^0\|_\XXX + \|\phi^0\|_{L^\infty} + 
  \lambda\|\partial_3 \phi^0\|_{L^\infty}^2 \,\le\, C\epsilon_2~,
$$
and so the smallness conditions on the perturbation in 
Theorem~\ref{asym_stab} imply those in Proposition~\ref{exproof}. 
We deduce that the solution of \reff{3DV2} with initial data
$\OO^0$ satisfies $\OO(x,t) = \OO^B(\xp;\rho+\phi(x_3,t)) +
\tilde\oomega(x,t)$ for some $\tilde\oomega \in \YYY_\mu(m)$, 
hence $\OO(x,t)$ converges exponentially in $X^2(m)^3$ toward 
the modulated vortex $\OO^B(\xp;\rho+\phi(x_3,t))$. On the
other hand, from \reff{phi_solve}, \reff{Gtdef}, we see that, 
for any $x_3 \in \real$, $\phi(x_3,t)$ converges toward the 
limiting value
$$
  \lim_{t \to \infty} \phi(x_3,t) \,=\, \frac{1}{\sqrt{2\pi}} 
  \int_{\real} e^{-\half y^2} \phi^0(y) \d y \,\equiv\, \delta \rho~,
$$
and that $\sup_{x_3 \in I}|\phi(x_3,t)-\delta\rho| = \cO(e^{-t})$ 
for any compact interval $I \subset \real$. Thus the difference 
$\|\OO^B(\cdot;\rho+\phi(x_3,t)) - \OO^B(\cdot;\rho+\delta\rho)
\|_{L^2(m)}$ will converge exponentially to zero as $t \to \infty$, 
uniformly for $x_3$ in any compact interval. Combining these
estimates, we obtain \reff{conv2}.


\section{Appendix}\label{appendix} 
 
 In this appendix we collect a number of technical estimates used in the main
 body of the paper.  They relate mostly to the behavior of the semigroup
 generated by the linearization of the vorticity equation around the Burgers
 vortex.  We also prove some estimates relating the vorticity field to the 
 corresponding velocity field defined by the Biot-Savart law.

\subsection{The one-dimensional Fokker-Planck operator}
\label{onedim}
 
Fix $\alpha > 0$, and consider the one-dimensional linear equation
\begin{equation}\label{Lalphadef}
  \partial_t \omega \,=\, \cL_{\alpha}\omega \,\equiv\,
  \partial_x^2 \omega + {\alpha \over 2}
  x \partial_x \omega + {\alpha \over 2}\omega~,
\end{equation}
where $x \in \real$ and $t \ge 0$. If $\omega(x,t) = \tilde\omega
(\sqrt{\alpha}x, \alpha t)$, then $\partial_t \tilde\omega
= \cL_1 \tilde\omega$ hence we could assume without loss of generality
that $\alpha = 1$. However for our purposes it is more convenient to
keep $\alpha > 0$ arbitrary.

The linear operator $\cL_\alpha$ is formally conjugated to the 
Hamiltonian of the harmonic operator in quantum mechanics:
$$
  e^{\frac{\alpha x^2}{8}} \,\cL_\alpha \,e^{-\frac{\alpha x^2}{8}}
  \,=\, L_\alpha \,\equiv\, \partial_x^2 -\frac{\alpha^2 x^2}{16}
  + \frac{\alpha}{4}~.
$$
As is well-known, the spectrum of $L_\alpha$ in $L^2(\real)$ is a 
sequence of simple eigenvalues:
$$
  \sigma(L_\alpha) \,=\, \Bigl\{-\frac{n\alpha}{2} \,\Big|\, 
  n = 0,1,2,\dots\Bigr\}~,
$$
and the associated eigenfunctions are the Hermite functions $h_n(x) =
e^{\alpha x^2/8} \partial_x^n e^{-\alpha x^2/4}$. This observation,
however, is not sufficient to determine the whole spectrum of
$\cL_\alpha$ because we want to consider this operator acting on a
space of functions with algebraic (rather than Gaussian) decay at
infinity.

For any $m \ge 0$ and $p \ge 1$ we define the space
$\LL^p(m) = \{f \in L^p(\real)\,|\, w^m f \in L^p(\real)\}$, where
$w(x) = (1+x^2)^{1/2}$. This Banach space is equipped with the 
natural norm
$$
   \|f\|_{\LL^p(m)} \,=\, \|w^m f\|_{L^p} \,=\, 
   \Bigl( \int_\real |w(x)^m f(x)|^p \d x\Bigr)^{1/p}~.
$$
The parameter $m$ determines the decay rate at infinity of functions
in $\LL^p(m)$. For instance, it is easy to verify that $\LL^2(m) 
\hookrightarrow L^1(\real)$ if (and only if) $m > 1/2$, because in that case
$w^{-m} \in L^2(\real)$ so that any $f \in \LL^2(m)$ satisfies
\begin{equation}\label{L1embedding}
   \int_\real |f(x)|\d x \,=\, \int_\real w(x)^m |f(x)| 
   w(x)^{-m}\d x \,\le\, \|w^m f\|_{L^2} \|w^{-m}\|_{L^2} \,=\, 
   C \|f\|_{\LL^2(m)}~,
\end{equation}
by H\"older's inequality. If $m > 1/2$, we thus define
$$
   \LL^2_0(m) \,=\, \Bigl\{f \in \LL^2(m) \,\Big|\, 
   \int_\real f(x)\d x = 0\Bigr\}~.
$$
This closed subspace of $\LL^2(m)$ is clearly invariant under the 
evolution defined by \reff{Lalphadef}.

In (\cite{gallay:2002a}, Appendix~A) it is shown that the spectrum 
of $\cL_\alpha$ in $\LL^2(m)$ is
\begin{equation}\label{sigmamdef}
  \sigma_m(\cL_\alpha) \,=\, \Bigl\{-\frac{n\alpha}{2} \,\Big|\, 
  n = 0,1,2,\dots\Bigr\} \,\cup\, \Bigl\{z \in \complex \,\Big|\, 
  \Re(z) \le \alpha\Bigl(\frac14 - \frac{m}{2}\Bigr)\Bigr\}~.
\end{equation}
Thus, in addition to the discrete spectrum of the harmonic 
oscillator, the operator $\cL_\alpha$ also has essential spectrum
due to the slow spatial decay of functions in $\LL^2(m)$. Note
however that this essential spectrum can be pushed far away 
from the imaginary axis by taking $m \ge 0$ sufficiently large.
Therefore, if $m$ is large, the relevant part of the spectrum 
of $\cL_\alpha$ is still given by the first few eigenvalues of 
$L_\alpha$. In particular zero is an isolated eigenvalue of 
$\cL_\alpha$ if $m > 1/2$, and the rest of the spectrum is 
strictly contained in the left-half plane. If we restrict 
ourselves to the invariant subspace $\LL^2_0(m)$, the
spectrum of $\cL_\alpha$ is unchanged except for the zero eigenvalue
(which is absent).

Equation \reff{Lalphadef} can be explicitly solved as
$\omega(t) = e^{t\cL_\alpha}\omega(0)$, where
$$
  \Bigl(e^{t\cL_\alpha}f\Bigr)(x) \,=\, \frac{e^{\alpha t/2}}
  {(4\pi a(t))^{1/2}} \int_\real e^{-\frac{(x-y)^2}{4a(t)}} 
  f(y e^{\alpha t/2})\d y~, \quad
  x \in \real~, \quad t > 0~,
$$
and $a(t) = (1 - e^{-\alpha t})/\alpha$. Using this expression, it 
is straightforward to verify that $e^{t\cL_\alpha}$ defines a 
strongly continuous semigroup in $\LL^2(m)$ for any $m \ge 0$. 
Moreover, $e^{t\cL_\alpha}$ maps $\LL^2_0(m)$ into $\LL^2_0(m)$ 
if $m > 1/2$, and the following estimates hold (see \cite{gallay:2002a},
Appendix~A): 

\begin{proposition}\label{1Dsemig}
If $m > 1/2$, the semigroup $e^{t\cL_\alpha}$ is uniformly bounded
in $\LL^2(m)$ for all $t \ge 0$. Moreover, if $m > 3/2$, there 
exists $C(m,\alpha) > 0$ such that, for all $f \in \LL^2_0(m)$, 
\begin{equation}\label{1Dest1}
  \|e^{t\cL_\alpha}f\|_{\LL^2(m)} \,\le\, C\,e^{-\alpha t/2} 
  \|f\|_{\LL^2(m)}~, \quad t \ge 0~.
\end{equation}
Finally, if $1 \le p \le 2$ and $m > 3/2$, then $e^{t\cL_\alpha}
\partial_x$ defines a bounded operator from $\LL^p(m)$ into
$\LL^2_0(m)$ and there exists $C(m,\alpha,p) > 0$ such that, 
for all $f \in \LL^p(m)$,
\begin{equation}\label{1Dest2}
  \|e^{t\cL_\alpha}\partial_x f\|_{\LL^2(m)} \,\le\, C
  \frac{e^{-\alpha t/2}}{a(t)^{\frac{1}{2p}+\frac14}} 
  \|f\|_{\LL^p(m)}~, \quad t > 0~,
\end{equation}
where $a(t) = (1 - e^{-\alpha t})/\alpha$. 
\end{proposition}

\subsection{Two-dimensional estimates}
\label{twodim}

We next consider the two-dimensional equation
\begin{equation}\label{Lalpha1alpha2}
  \partial_t \omega \,=\,  \cL_{\alpha_1,\alpha_2}\omega \,\equiv\,
  \Delta \omega + {\alpha_1 \over 2} x_1 \partial_1 \omega 
  + {\alpha_2 \over 2} x_2 \partial_2 \omega + 
  \frac{\alpha_1 + \alpha_2}{2}\omega~,
\end{equation}
where $x \in \real^2$, $ t\ge 0$, and $\alpha_1 \ge \alpha_2 > 0$.
In the particular case where $\alpha_1 = 1+\lambda$ and $\alpha_2 
= 1 - \lambda$ for some $\lambda \in [0,1)$, we see that 
$\cL_{\alpha_1,\alpha_2} = \cLp + \lambda \cM$, where 
$\cLp, \cM$ are defined in \reff{LMdef}. Note that the parameters
$\alpha_1,\alpha_2$ cannot be eliminated by a rescaling, unless 
$\alpha_1 = \alpha_2$. 

We study the operator $\cL_{\alpha_1,\alpha_2}$ in the weighted space
\begin{equation}\label{Lpmdef}
  L^p(m) \,=\, \{f \in L^p(\real^2) \,|\, b^m f \in L^p(\real^2)\}~,
  \quad \|f\|_{L^p(m)} = \|b^m f\|_{L^p}~,
\end{equation}
where $p \ge 1$, $m \ge 0$, and $b(x_1,x_2) = w(x_1)w(x_2) = 
(1+x_1^2)^{1/2}(1+x_2^2)^{1/2}$. It is clear that $L^2(m) = \LL^2(m) 
\otimes \LL^2(m)$, where $\LL^2(m)$ is the one-dimensional space
defined in the previous paragraph and $\otimes$ denotes the tensor
product of Hilbert spaces, see \cite{reed:1972}. Comparing the 
definitions \reff{Lalphadef}, \reff{Lalpha1alpha2}, we see that 
our operator can be decomposed as $\cL_{\alpha_1,\alpha_2} = 
\cL_{\alpha_1} \otimes \oone + \oone\otimes \cL_{\alpha_2}$, where
$\oone$ denotes the identity operator. It follows that the spectrum of 
$\cL_{\alpha_1,\alpha_2}$ in $L^2(m)$ is just the sum
$$
  \sigma_m(\cL_{\alpha_1,\alpha_2}) \,=\, \sigma_m(\cL_{\alpha_1}) +  
  \sigma_m(\cL_{\alpha_2})~,
$$
where $\sigma_m(\cL_{\alpha_i})$ is given by \reff{sigmamdef} for
$i = 1,2$. In particular, zero is an isolated eigenvalue of 
$\cL_{\alpha_1,\alpha_2}$ if $m > 1/2$, and if $m > 3/2$ there 
exists $\mu > \alpha_2/2$ such that
$$
  \sigma_m(\cL_{\alpha_1,\alpha_2}) \,\subset\, 
  \Bigl\{0 \,,\, -\frac{\alpha_2}{2}\Bigr\} \,\cup\, 
  \Bigl\{z \in \complex \,\Big|\, \Re(z) \le -\mu\Bigr\}~.
$$
(Recall that we assumed $\alpha_1 \ge \alpha_2$.) Moreover, if 
$m > 1/2$, the subspace $L^2_0(m)$ defined by \reff{L20mdef}
is invariant under the action of $\cL_{\alpha_1,\alpha_2}$, and
the restriction of $\cL_{\alpha_1,\alpha_2}$ to $L^2_0(m)$ has
spectrum $\sigma_m(\cL_{\alpha_1,\alpha_2}) \setminus \{0\}$. 
Thus $\cL_{\alpha_1,\alpha_2}$ is invertible in $L^2_0(m)$ if 
$m > 1/2$, with bounded inverse. 

The semigroup generated by $\cL_{\alpha_1,\alpha_2}$ satisfies 
$e^{t\cL_{\alpha_1,\alpha_2}} = e^{t\cL_{\alpha_1}} \otimes 
e^{t\cL_{\alpha_2}}$. Thus, using Proposition~\ref{1Dsemig}, we 
immediately obtain the following estimates:

\begin{proposition}\label{2Dsemig}
If $m > 1/2$, the semigroup $e^{t\cL_{\alpha_1,\alpha_2}}$ is
uniformly bounded in $L^2(m)$ for all $t \ge 0$. Moreover, 
if $m > 3/2$, there exists $C(m,\alpha_1,\alpha_2) > 0$ such that, 
for all $f \in L^2_0(m)$, 
\begin{equation}\label{2Dest0}
  \|e^{t\cL_{\alpha_1,\alpha_2}}f\|_{L^2(m)} \,\le\, C\,e^{-\alpha_2 t/2} 
  \|f\|_{L^2(m)}~, \quad t \ge 0~.
\end{equation}
Finally, if $1 \le p \le 2$ and $m > 3/2$, then 
$e^{t\cL_{\alpha_1,\alpha_2}}\partial_k$ defines a bounded operator
from $L^p(m)$ into $L^2_0(m)$ for $k = 1,2$, and there exists 
$C(m,\alpha_1,\alpha_2,p) > 0$ such that, for all $f \in L^p(m)$,
\begin{eqnarray}\label{2Dest1}
  \|e^{t\cL_{\alpha_1,\alpha_2}}\partial_1 f\|_{L^2(m)} \,\le\, C
    \frac{e^{-\alpha_1 t/2}}{a_1(t)^{\frac{1}{2p}+\frac14}
    a_2(t)^{\frac{1}{2p}-\frac14}}\|f\|_{L^p(m)}~, \quad t > 0~,
   \\ \label{2Dest2}
  \|e^{t\cL_{\alpha_1,\alpha_2}}\partial_2 f\|_{L^2(m)} \,\le\, C
    \frac{e^{-\alpha_2 t/2}}{a_1(t)^{\frac{1}{2p}-\frac14}
    a_2(t)^{\frac{1}{2p}+\frac14}}\|f\|_{L^p(m)}~, \quad t > 0~,
\end{eqnarray}
where 
$$
   a_i(t) \,=\, \frac{1 - e^{-\alpha_i t}}{\alpha_i} \,=\, 
   \int_0^t e^{-\alpha_i s}\d s~, \quad i = 1,2~.
$$
\end{proposition}

\begin{remark}\label{noderiv}
For $p \in [1,2]$ and $m > 1/2$, we also have the following bound:
$$
  \|e^{t\cL_{\alpha_1,\alpha_2}} f\|_{L^2(m)} \,\le\, 
    \frac{C}{a_1(t)^{\frac{1}{2p}-\frac14} a_2(t)^{\frac{1}{2p}-\frac14}}
  \|f\|_{L^p(m)}~, \quad t > 0~.
$$
\end{remark}

\medskip
We conclude this paragraph with a short discussion of the 
two-dimensional Biot-Savart law:
\begin{equation}\label{BS2Dbis}
  \uu(x) \,=\, \frac{1}{2\pi} \inttwo
  \frac{1}{|x-y|^2}\pmatrix{y_2 - x_2 \cr x_1 - y_1}
  \omega(y)\d y~, \quad x \in \real^2~.
\end{equation}

\begin{proposition}\label{BS2Destimates}
Let $\uu$ be the velocity field defined from $\omega$ via the 
Biot-Savart law \reff{BS2Dbis}.\\[1mm]
i) If $\omega \in L^p(\real^2)$ for some $p \in (1,2)$, then
$\uu \in L^q(\real^2)$ where $\frac1q = \frac1p - \frac12$, and there 
exists $C(p) > 0$ such that $\|\uu\|_{L^q} \le C\|\omega\|_{L^p}$.\\[1mm]
ii) If $\omega \in L^p(\real^2) \cap L^q(\real^2)$ for some $p \in 
[1,2)$ and some $q \in (2,+\infty]$ then $\uu \in C^0_b(\real^2)$ 
and there exists $C(p,q) > 0$ such that
$$
  \|\uu\|_{L^{\infty}} \,\le\, C \|\omega\|_{L^p}^\alpha 
  \|\omega\|_{L^q}^{1-\alpha}~, \quad \hbox{where}\quad 
  \frac12 \,=\, \frac{\alpha}p + \frac{1-\alpha}q~. 
$$
\end{proposition}

\proof Assertion i) is a direct consequence of the 
Hardy-Littlewood-Sobolev inequality \cite{lieb:1997}. For a proof
of ii), see for instance (\cite{gallay:2002a}, Lemma~2.1). 
\QED

\medskip
We deduce from Proposition~\ref{BS2Destimates} the following 
useful bound on the product $\uu \omega$:

\begin{corollary}\label{2Dprod}
Assume that $\omega_1, \omega_2 \in L^2(m)$ for some $m > 1/2$, 
and let $\uu_1$ be the velocity field obtained from $\omega_1$
via the Biot-Savart law \reff{BS2Dbis}. Then $\uu_1 \omega_2 \in 
L^p(m)$ for all $p \in (1,2)$, and there exists $C(m,p) > 0$ 
such that 
$$
  \|\uu_1 \omega_2\|_{L^p(m)} \le C \|\omega_1\|_{L^2(m)}
  \|\omega_2\|_{L^2(m)}~.
$$
\end{corollary}

\proof 
Assume that $1 < p < 2$. By H\"older's inequality
$$
  \|\uu_1 \omega_2\|_{L^p(m)} \,=\, \|b^m \uu_1 \omega_2
  \|_{L^p}  \,\le\, \|\uu_1\|_{L^q} \|b^m \omega_2\|_{L^2}~,
  \quad {\rm where} \quad \frac1q \,=\, \frac1p - \frac12~.
$$
Now $\|\uu_1\|_{L^q} \le C \|\omega_1\|_{L^p}$ by
Proposition~\ref{BS2Destimates}, and $\|\omega_1\|_{L^p} 
\le C \|\omega_1\|_{L^2(m)}$ because $L^2(m) \hookrightarrow 
L^p(\real^2)$ for $p \in [1,2]$ if $m > 1/2$. This gives the
desired result. \QED

\subsection{The three-dimensional semigroup}
\label{threedim}

This section is devoted to the three-dimensional equation 
\begin{equation}\label{Veq}
  \partial_t \omega \,=\,  \hat\cL_{\alpha_1,\alpha_2}\omega \,\equiv\,
  \Delta \omega + {\alpha_1 \over 2}
  x_1 \partial_1 \omega + {\alpha_2 \over 2}
  x_2 \partial_2 \omega - x_3 \partial_3 \omega + 
  \frac{\alpha_1+\alpha_2}{2}\omega~,
\end{equation}
where $x \in \real^3$, $t \ge 0$, and $\alpha_1 \ge \alpha_2 > 0$.
In the particular case where $\alpha_1 = 1+\lambda$ and $\alpha_2 = 
1-\lambda$, we have $\hat\cL_{\alpha_1,\alpha_2} = \cL + 1$ where 
$\cL$ is defined in \reff{cLdef}. 

It is important to realize that the evolution defined by \reff{Veq}
is essentially contracting in the transverse variables $\xp = 
(x_1,x_2)$ and expanding in the axial variable $x_3$. This is due
to the signs of the advection terms, which in turn originate in 
our choice of the straining flow \reff{strain}. For this reason 
we can assume that the solutions of \reff{Veq} decay to zero as
$|\xp| \to \infty$, but we cannot impose any decay in the $x_3$
variable (otherwise the solutions will not stay uniformly bounded 
for all times in the corresponding norm). This motivates the following
choice of our function space. For $p \ge 1$ and $m \ge 0$, we
introduce the Banach space
\begin{equation}\label{Xpmdef}
  X^p(m) \,\equiv\, C^0_b(\real,L^p(m)) \,=\, \{\omega : \real 
  \to L^p(m)\,|\, \omega \hbox{ is bounded and continuous}\}
\end{equation}
equipped with the norm
$$
  \|\omega\|_{X^p(m)} \,=\, \sup_{x_3\in\real} 
  \|\omega(x_3)\|_{L^p(m)}~.
$$
For any $n \in \intplus^*$ we also define the seminorm
\begin{equation}\label{pseminorm}
  |\omega|_{X^p_n(m)} \,=\, \sup_{|x_3|\le n} \|\omega(x_3)\|_{L^p(m)}~,
\end{equation}
and we denote by $X^p_\loc(m)$ the space $X^p(m)$ equipped with 
the topology defined by the family of seminorms \reff{pseminorm}
for $n \in \intplus^*$. For later use, we observe that the ball 
$\{\omega \in X^p(m)\,|\, \|\omega\|_{X^p(m)} \le R\}$ is closed in 
$X^p_\loc(m)$ for any $R > 0$. 

At least formally, the space $X^2(m)$ can be thought of as the 
tensor product $C^0_b(\real) \otimes L^2(m)$, i.e. the space generated 
by linear combinations of elements of the form $\omega(\xp,x_3) 
= f(x_3)g(\xp)$, with $f \in C^0_b(\real)$ and $g \in L^2(m)$. 
In this picture, the linear operator defined by \reff{Veq} can
be decomposed as $\hat\cL_{\alpha_1,\alpha_2} = \hat\cL_3 \otimes 
\oone + \oone \otimes \cL_{\alpha_1,\alpha_2}$, where 
$\cL_{\alpha_1,\alpha_2}$ is defined in \reff{Lalpha1alpha2} and 
$\hat\cL_3$ is the one-dimensional operator $\hat\cL_3 f =
\partial_3^2 f -x_3\partial_3 f$. It is readily verified
that $\hat\cL_3$ generates a semigroup in $C^0_b(\real)$ given by 
the explicit formula
\begin{equation}\label{L3exp}
 (e^{t\hat\cL_3}f)(x_3) \,=\, (G_t * f)(x_3 e^{-t})~, \quad x_3 \in
 \real~, \quad t > 0~,
 \end{equation}
where $G_t$ is defined in \reff{Gtdef}, and we know from 
Section~\ref{twodim} that $\cL_{\alpha_1,\alpha_2}$ generates a 
strongly continuous semigroup in $L^2(m)$. Thus we expect that 
$\hat\cL_{\alpha_1,\alpha_2}$ will generate a semigroup 
$\{S_t\}_{t\ge 0}$ in $X^2(m)$ given by $S_t = e^{t\hat\cL_3} \otimes 
e^{t\cL_{\alpha_1,\alpha_2}}$, or explicitly
\begin{equation}\label{sg1}
  (S_t \omega)(x_3) \,=\, \int_\real G_t(x_3 e^{-t} - y_3) 
  \Bigl(e^{t\cL_{\alpha_1,\alpha_2}}\omega(y_3)\Bigr)\d y_3~, 
  \quad x_3 \in \real~, \quad t > 0~.
\end{equation}
We shall prove that these heuristic considerations are indeed 
correct in the sense that \reff{sg1} defines a semigroup of bounded
operators in $X^2(m)$ with the property that $\omega(t) = S_t \omega$
is the solution of \reff{Veq} with initial data $\omega \in X^2(m)$.
However, the map $t \mapsto S_t \omega$ is not continuous in the 
topology of $X^2(m)$, but only in the (weaker) topology of 
$X^2_\loc(m)$. This is due to the fact that equation \reff{Veq} has 
``infinite speed of propagation'' in the sense that the advection term 
in the vertical variable is unbounded, see Remark~\ref{notcontinuous}.

\begin{proposition}\label{3Dsemig} 
For any $m \ge 0$, the family $\{S_t\}_{t \ge 0}$ defined by 
\reff{sg1} and $S_0 = \oone$ is a semigroup of bounded linear 
operators on $X^2(m)$. If $\omega_0 \in X^2(m)$ and $\omega(t) = 
S_t \omega_0$, then $\omega : [0,+\infty) \to X^2_\loc(m)$ is 
continuous, and $\omega(t)$ solves \reff{Veq} for $t > 0$. 
For any $R > 0$, if $B_R \,=\, \{f \in X^2(m)\,|\, \|f\|_{X^2(m)} 
\le R\}$ is equipped with the topology of $X^2_\loc(m)$, then 
$S_t : B_R \to X^2_\loc(m)$ is continuous, uniformly in time on compact 
intervals. Moreover:\\[1mm]
i) If $m > 1/2$ then $S_t$ is uniformly bounded on $X^2(m)$ for all 
$t \ge 0$.\\[1mm]
ii) If $m > 3/2$ there exists $C(m,\alpha_1,\alpha_2) > 0$ such that, 
for all $\omega$ in the subspace $X^2_0(m)$ defined in \reff{X20mdef},
\begin{eqnarray}\label{3Dest0}
  \|S_t \omega\|_{X^2(m)} \,\le\, C \,e^{-\frac{\alpha_2}{2} t} 
  \|\omega\|_{X^2(m)}~, \quad t \ge 0~.
\end{eqnarray}
iii) If $p \in [1,2]$ and $m > 3/2$, $S_t\partial_k$ defines a bounded
operator from $X^p(m)$ into $X^2(m)$ for $t > 0$ and $k = 1,2,3$, 
and there exists $C(m,\alpha_1,\alpha_2,p) > 0$ such that 
\begin{eqnarray} \label{3Dest1}
  \|S_t \partial_1 \omega\|_{X^2(m)} 
  &\le& C \,\frac{e^{-\alpha_1 t/2}}{a_1(t)^{\frac{1}{2p}
  +\frac14} a_2(t)^{\frac{1}{2p}-\frac14}} \,\|\omega\|_{X^p(m)}~,
  \\ \label{3Dest2}
  \|S_t \partial_2 \omega\|_{X^2(m)}
  &\le& C \,\frac{e^{-\alpha_2 t/2}}{a_1(t)^{\frac{1}{2p}
  -\frac14} a_2(t)^{\frac{1}{2p}+\frac14}} \,\|\omega\|_{X^p(m)}~,
  \\ \label{3Dest3}
  \|S_t \partial_3 \omega\|_{X^2(m)}
  &\le& \frac{C}{\sqrt{1-e^{-2t}}} 
  \,\frac{1}{a_1(t)^{\frac{1}{2p}-\frac14} 
  a_2(t)^{\frac{1}{2p}-\frac14}} \,\|\omega\|_{X^p(m)}~,
\end{eqnarray}
where $a_1(t), a_2(t)$ are as in Proposition~\ref{2Dsemig}. 
\end{proposition}

\proof
We first rewrite \reff{sg1} in a slightly more convenient form. 
By \reff{Gtdef} we have $G_t(y) = c(t)^{-1/2} G(c(t)^{-1/2}y)$, where
$c(t) = 1-e^{-2t}$ and $G(z) = (2\pi)^{-1/2}e^{-z^2/2}$. Thus setting
$y_3 = x_3 e^{-t} + c(t)^{1/2}z_3$ in \reff{sg1}, we obtain the 
equivalent formula
\begin{equation}\label{sg2}
  (S_t \omega)(x_3) \,=\, \int_\real G(z_3) \Bigl(e^{t\cL_{\alpha_1,
  \alpha_2}}\omega(x_3 e^{-t} + c(t)^{1/2}z_3)\Bigr)\d z_3~, 
  \quad x_3 \in \real~, \quad t \ge 0~.
\end{equation}
Fix $m \ge 0$. If $\omega \in X^2(m)$, then for any $t \ge 0$ the 
map $x_3 \mapsto e^{t\cL_{\alpha_1,\alpha_2}}\omega(x_3)$ also belongs
to $X^2(m)$, because $e^{t\cL_{\alpha_1,\alpha_2}}$ is a bounded
operator on $L^2(m)$ by Proposition~\ref{2Dsemig}. Thus it follows 
immediately from \reff{sg2} that $S_t \omega \in X^2(m)$ and
\begin{equation}\label{Stbd1}
  \|S_t \omega\|_{X^2(m)} \,\le\, \sup_{x_3 \in \real} 
  \|e^{t\cL_{\alpha_1,\alpha_2}}\omega(x_3)\|_{L^2(m)} \,\le\, 
  N_m(t) \|\omega\|_{X^2(m)}~,
\end{equation}
where $N_m(t) = \|e^{t\cL_{\alpha_1,\alpha_2}}\|_{L^2(m)\to L^2(m)}$. 
The semigroup formula $S_{t_1+t_2} = S_{t_1}S_{t_2}$ is easily verified 
using \reff{sg1}, Fubini's theorem, and the identity
$$
  \int_{\real} G_{t_1}(xe^{-t_1}-y) G_{t_2}(ye^{-t_2}-z)\d y
  \,=\, G_{t_1+t_2}(xe^{-(t_1+t_2)}-z)~.
$$
Thus $\{S_t\}_{t \ge 0}$ is a semigroup of bounded operators in 
$X^2(m)$. 

On the other hand, by \reff{sg2}, we have for all $k \in \intplus$:
$$
  \|S_t \omega(x_3)\|_{L^2(m)} \,\le\, \int_{-k}^k G(z_3)N_m(t)
  \|\omega(x_3 e^{-t} + c(t)^{1/2}z_3)\|_{L^2(m)}\d z_3 + 
  N_m(t) \|\omega\|_{X^2(m)} \epsilon_k~,
$$
where $\epsilon_k = \int_{|z|\ge k}G(z)\d z \to 0$ as $k \to \infty$. 
Since $|x_3 e^{-t} + c(t)^{1/2}z_3| \le n+k$ whenever $|x_3| \le n$ 
and $|z_3| \le k$, we deduce that for all $n \in \intplus^*$:
\begin{equation}\label{Stbd2}
  |S_t \omega|_{X^2_n(m)} \,\le\, N_m(t)\Bigl(|\omega|_{X^2_{n+k}(m)}
  + \epsilon_k \|\omega\|_{X^2(m)}\Bigr)~.
\end{equation}
This bound implies that $S_t : B_R \to X^2_\loc(m)$ is continuous, 
uniformly in time on compact intervals. 

Furthermore, if $\omega \in X^2(m)$ and $t > 0$, we have
\begin{eqnarray*}
  (S_t \omega - \omega)(x_3) &=& \int_\real G(z_3)
  e^{t\cL_{\alpha_1,\alpha_2}}\Bigl(\omega(x_3 e^{-t} + c(t)^{1/2}z_3)
  -\omega(x_3)\Bigr)\d z_3 \\
  &+& e^{t\cL_{\alpha_1,\alpha_2}}\omega(x_3)-\omega(x_3)~,
\end{eqnarray*}
hence proceeding as above we find for all $n,k \in \intplus^*$:
\begin{eqnarray*}
  |S_t \omega - \omega|_{X^2_n(m)} &\le& \int_{-k}^k G(z_3)
  N_m(t) \sup_{|x_3|\le n} \|\omega(x_3 e^{-t} + c(t)^{1/2}z_3)
  -\omega(x_3)\|_{L^2(m)} \d z_3 \\
  &+& 2 N_m(t) \epsilon_k \|\omega\|_{X^2(m)} +  
  \sup_{|x_3|\le n} \|e^{t\cL_{\alpha_1,\alpha_2}}\omega(x_3)
  -\omega(x_3)\|_{L^2(m)}~.
\end{eqnarray*}
The last term goes to zero as $t \to 0+$ because $e^{t\cL_{\alpha_1,
\alpha_2}}$ is a strongly continuous semigroup on $L^2(m)$ and 
$\omega : [-n,n] \to L^2(m)$ is continuous (hence has compact 
range). Similarly, for each $k \in \intplus$, the integral goes to
zero as $t \to 0+$ by Lebesgue's dominated convergence theorem, 
because $\omega : [-n{-}k,n{+}k] \to L^2(m)$ is uniformly continuous. 
Since $\epsilon_k \to 0$ as $k \to \infty$, we conclude that 
$S_t\omega \to \omega$ in $X^2_\loc(m)$ as $t \to 0+$. Then, using 
the semigroup property, we deduce that the map $t \mapsto S_t \omega$ 
is continuous to the right at any $t \ge 0$. Finally, if $t >
\epsilon > 0$,  we have $S_t\omega - S_{t-\epsilon}\omega = 
S_{t-\epsilon}(S_\epsilon \omega - \omega)$, hence by \reff{Stbd2}
$$
  |S_t \omega - S_{t-\epsilon}\omega|_{X^2_n(m)} \,\le\, 
  N_m(t-\epsilon)\Bigl(|S_\epsilon\omega - \omega|_{X^2_{n+k}(m)} + 
  \epsilon_k (1+N_m(\epsilon))\|\omega\|_{X^2(m)}\Bigr)~.
$$
Since $|S_\epsilon\omega - \omega|_{X^2_{n+k}(m)} \to 0$ as 
$\epsilon \to 0+$ for all $k,n \in \intplus^*$, this shows that
$t \mapsto S_t \omega$ is also continuous to the left at any 
$t > 0$. 

Next, using \reff{sg2}, \reff{L3exp}, and the explicit formula
for $e^{t\cL_{\alpha_1,\alpha_2}}$, it is rather straightforward to verify 
that, for any $\omega \in X^2(m)$, the map $(\xp,x_3,t) \mapsto 
\omega(\xp,x_3,t) = ((S_t \omega)(x_3))(\xp)$ is smooth and satisfies
\reff{Veq} for $t > 0$. Thus $S_t \omega$ is indeed the solution
of \reff{Veq} with initial data $\omega$. 

\smallskip
It remains to establish the decay properties of $S_t$:\\[1mm]
i) If $m > 1/2$, we know from Proposition~\ref{2Dsemig} that
$N_m(t) \le C$ for all $t \ge 0$, hence $\{S_t\}_{t \ge 0}$
is uniformly bounded on $X^2(m)$ by \reff{Stbd1}.\\[1mm]
ii) If $m > 3/2$ and $\omega \in X^2_0(m)$, then $\omega(x_3) \in 
L^2_0(m)$ for all $x_3 \in \real$ and \reff{3Dest0} follows
immediately from \reff{Stbd1} and \reff{2Dest0}.\\[1mm]
iii) If $k = 1,2$, we define $S_t\partial_k$ by 
\begin{equation}\label{sg3}
  (S_t \partial_k \omega)(x_3) \,=\, \int_\real G(z_3) 
  \Bigl(e^{t\cL_{\alpha_1,\alpha_2}}\partial_k \omega(x_3 e^{-t} 
  + c(t)^{1/2}z_3)\Bigr)\d z_3~. 
\end{equation}
If $m > 3/2$ and $p \in [1,2]$, we know from Proposition~\ref{2Dsemig}
that $e^{t\cL_{\alpha_1,\alpha_2}}\partial_k$ is a bounded operator 
from $L^p(m)$ into $L^2_0(m)$ satisfying \reff{2Dest1} or
\reff{2Dest2}. Thus the formula \reff{sg3} defines a bounded operator
from $X^p(m)$ into $X^2_0(m)$ and
$$
  \|S_t \partial_k \omega\|_{X^2(m)} \,\le\, \sup_{x_3 \in \real}
  \|\e^{t\cL_{\alpha_1,\alpha_2}} \partial_k \omega(x_3)\|_{L^2(m)}~.
$$
Thus \reff{3Dest1}, \reff{3Dest2} follow immediately from 
\reff{2Dest1}, \reff{2Dest2}. Finally we define $S_t\partial_3$ by 
\begin{equation}\label{sg4}
  (S_t \partial_3 \omega)(x_3) \,=\, -\frac{1}{c(t)^{1/2}}
  \int_\real \partial_3 G(z_3) \Bigl(e^{t\cL_{\alpha_1,
  \alpha_2}}\omega(x_3 e^{-t} + c(t)^{1/2}z_3)\Bigr)\d z_3~.
\end{equation}
We know that $e^{t\cL_{\alpha_1,\alpha_2}}$ is a bounded operator from
$L^p(m)$ into $L^2(m)$, see Remark~\ref{noderiv}. Since $\partial_3 G
\in L^1(\real^2)$, we thus find
\begin{eqnarray*}
  \|S_t \partial_3 \omega\|_{X^2(m)} &\le& \frac{C}{c(t)^{1/2}}
  \sup_{x_3 \in \real} \|e^{t\cL_{\alpha_1,\alpha_2}}\omega(x_3)
  \|_{L^2(m)}\\  &\le& \frac{C}{c(t)^{1/2}} \frac{1}
  {a_1(t)^{\frac{1}{2p}-\frac14} a_2(t)^{\frac{1}{2p}-\frac14}}
  \|\omega\|_{X^p(m)}~,
\end{eqnarray*}
which is \reff{3Dest3}. This concludes the proof. \QED

\begin{corollary}\label{3Dint}
Let $T > 0$ and let $f : [0,T] \to X^2(m)$ be a bounded function 
satisfying $f \in C^0([0,T],X^2_\loc(m))$. Then the map $F : 
[0,T] \to X^2(m)$ defined by
$$
  F(t) \,=\, \int_0^t S_{t-s}f(s)\d s~, \quad 0 \le t \le T~,
$$
satisfies $F \in C^0([0,T],X^2_\loc(m))$, and $\|F(t)\|_{X^2(m)} 
\le \int_0^t N_m(t-s)\|f(s)\|_{X^2(m)}\d s$, where $N_m$ is as 
in \reff{Stbd1}. 
\end{corollary}

\proof For any $t \in [0,T]$, we define $\psi_t : [0,t] \to X^2(m)$ 
by $\psi_t(s) \,=\, S_{t-s}f(s)$. Since $f \in C^0([0,T],X^2_\loc(m))$ 
and since the semigroup $S_t$ is continuous on $X^2_\loc(m)$ as 
described in Proposition~\ref{3Dsemig}, it is easy to verify that 
the map $\psi_t : [0,t] \to X^2_\loc(m)$ is also continuous. 
As $X^2_\loc(m)$ is a subspace of the Fr\'echet space $C^0(\real,L^2(m))$, 
the integral $\int_0^t \psi_t(s)\d s$ can be defined as in 
(\cite{rudin:1991}, Theorem~3.17). However, in the present case, we 
can also use the following ``pedestrian'' construction (which agrees 
with the general one). For any $n \in \intplus^*$, we define
$$
  F_n(t) \,=\, \int_0^t \chi_n \psi_t(s)\d s \,=\, 
  \int_0^t \chi_n S_{t-s}f(s)\d s~,
$$
where $\chi_n$ denotes the map $x_3 \mapsto \oone_{[-n,n]}(x_3)$. 
Clearly $\chi_n \psi_t(s)$ is a continuous function of $s \in [0,t]$
with values in the Banach space $C^0([-n,n],L^2(m))$, hence $F_n(t) 
\in C^0([-n,n],L^2(m))$ can be defined for any $t \in [0,T]$ as 
a Banach-valued Riemann integral. Using again the continuity of the 
semigroup $S_t$ one finds that $F_n : [0,T] \to C^0([-n,n],L^2(m))$
is continuous and satisfies
$$
  |F_n(t)|_{X^2_n(m)} \,\le\, \int_0^t N_m(t-s)\|f(s)\|_{X^2(m)} 
  \d s \,\le\, C(T) \int_0^t \|f(s)\|_{X^2(m)} \d s~.
$$
(Note that $t \mapsto N_m(t)$ and $t \mapsto \|f(t)\|_{X^2(m)}$ are
lower semicontinuous, hence measurable.) Now, for each $t \in [0,T]$,
it is clear that $(F_m(t))(x_3) = (F_n(t))(x_3)$ if $|x_3| \le n \le
m$, hence there is a unique $F(t) \in C^0(\real,L^2(m))$ such that
$(F(t))(x_3) = (F_n(t))(x_3)$ whenever $|x_3| \le n$. By construction,
$\|F(t)\|_{X^2(m)} \le \int_0^t N_m(t-s)\|f(s)\|_{X^2(m)} \d s$ for
all $t \in [0,T]$, and $F \in C^0([0,T],X^2_\loc(m))$. \QED

\begin{remark}\label{3Dint2}
Similarly, if $p \in (1,2]$ and $k \in \{1,2,3\}$, 
Proposition~\ref{3Dsemig} implies that, if $f : [0,T] \to X^p(m)$ 
is bounded in $X^p(m)$ and continuous in $X^p_\loc(m)$, the
map $F :  [0,T] \to X^2(m)$ defined by
$$
  F(t) \,=\, \int_0^t S_{t-s}\partial_k f(s)\d s~, 
  \quad 0 \le t \le T~,
$$
is bounded in $X^2(m)$ and continuous in $X^2_\loc(m)$. In that
case, for each $n \in \intplus^*$ and each $t \in (0,T]$, 
$F_n(t) \equiv \chi_n F(t)$ is defined by a ``generalized''
Riemann integral, because the integrand has a singularity 
at $s = t$.  
\end{remark}

\subsection{The three-dimensional Biot-Savart law}
\label{biotsavart}

In this final section, we discuss the three-dimensional Biot-Savart 
law, namely
\begin{equation}\label{3DBS}
   \uu(x) \,=\, -\frac{1}{4\pi} \int_{\real^3} 
   \frac{(\xx-\yy) \times \oomega(y)}{|x-y|^3}\d y~, \quad
   x \in \real^3~.
\end{equation}
We first prove the analogue of Proposition~\ref{BS2Destimates} 
in the spaces $X^p(m)$ defined by \reff{Xpmdef}. 

\begin{proposition}\label{velocity-vorticity}
Fix $m > 1/2$. If $\oomega \in X^2(m)$, the velocity field given 
by \reff{3DBS} satisfies $\uu \in X^q(0)$ for all $q \in (2,\infty)$, 
and there exists $C(m,q) > 0$ such that $\|\uu\|_{X^q(0)} \le 
C\|\oomega\|_{X^2(m)}$. 
\end{proposition}

\proof
Assume that $\omega \in X^2(m)$ for some $m > 1/2$. For all $x = 
(\xp,x_3) \in \real^3$, we have by Fubini's theorem
$$
  |\uu(\xp,x_3)| \,\le\, C \int_{\real^3} \frac{|\oomega(\yp,y_3)|}
  {|\xp-\yp|^2 + (x_3-y_3)^2} \d\yp \d y_3 \,=\,
  C \int_{\real} F(\xp;x_3,y_3)\d y_3~,
$$
where
$$
  F(\xp;x_3,y_3) \,=\, \inttwo\frac{|\oomega(\yp,y_3)|}
  {|\xp-\yp|^2 + (x_3-y_3)^2} \d\yp~.
$$
By Minkowski's inequality, it follows that
\begin{equation}\label{uubd1}
  \|\uu(\cdot,x_3)\|_{L^q(\real^2)} \,\le\, C \int_{\real}  
  \|F(\cdot;x_3,y_3)\|_{L^q(\real^2)} \d y_3~.
\end{equation}
If $2 < q < \infty$, we shall show that there exists $H_q \in
L^1(\real)$ and $C > 0$ such that
\begin{equation}\label{Heq}
  \|F(\cdot;x_3,y_3)\|_{L^q(\real^2)} \,\le\, C 
  \|\oomega(\cdot,y_3)\|_{L^2(m)} H_q(x_3-y_3)~, \quad x_3,y_3 \in \real~.
\end{equation}
Together with \reff{uubd1}, this gives $\|\uu(\cdot,x_3)\|_{L^q(\real^2)}
\le C \|\oomega\|_{X^2(m)}$ for all $x_3 \in \real$, which is the 
desired bound. Since the Biot-Savart law is invariant under spatial 
translations, the same arguments show that, for all $x_3 \in \real$, 
$$
  \|\uu(\cdot,x_3+\epsilon) - \uu(\cdot,x_3)\|_{L^q(\real^2)}
  \,\le\, C \int_\real \|\oomega(\cdot,y_3+\epsilon)- 
  \oomega(\cdot,y_3)\|_{L^2(m)} H_q(x_3-y_3)\d y_3~.
$$
As $\epsilon \to 0$ the right-hand side converges to zero by 
Lebesgue's dominated convergence theorem, thus $\uu \in C^0_b(\real,
L^q(\real^2)) \equiv X^q(0)$. 

To prove \reff{Heq}, for any $a \in \real$ we define $f_a(\yp) = 
(|\yp|^2 + a^2)^{-1}$. If $a \neq 0$, then $f_a \in L^r(\real^2)$ for 
all $r > 1$, and there exists $C_r > 0$ such that 
$$
   \|f_a\|_{L^r(\real^2)} \,\le\, \frac{C_r}{|a|^{2-\frac2r}}~.
$$ 
Since $F(\cdot;x_3,y_3) = |\oomega(\cdot,y_3)| * f_{x_3-y_3}$, Young's
inequality implies
\begin{equation}\label{Hbdd1}
  \|F(\cdot;x_3,y_3)\|_{L^q(\real^2)} \,\le\, C \|\oomega(\cdot,y_3)
  \|_{L^2(\real^2)} \|f_{x_3-y_3}\|_{L^p(\real^2)} \,\le\, 
  \frac{C \|\oomega(\cdot,y_3)\|_{L^2(\real^2)}}{|x_3-y_3|^{2-\frac2p}}~,
\end{equation}
where $1 + \frac1q = \frac12 + \frac1p$, hence $2 -\frac2p = 1 -
\frac2q < 1$. On the other hand, by H\"older's inequality,
\begin{eqnarray*}
  F(\xp;x_3,y_3) &=& \inttwo b^m(\yp) |\oomega(\yp,y_3)|
  \,\frac{1}{b^m(\yp) (|\xp{-}\yp|^2 + (x_3{-}y_3)^2)} \d\yp \\
  &\le& \|b^m\oomega(\cdot,y_3)\|_{L^2(\real^2)} 
  \Bigl(\frac{1}{b^{2m}} * f^2_{x_3-y_3} \Bigr)^{1/2}(\xp)~.
\end{eqnarray*}
Using Young's inequality again, we obtain
\begin{eqnarray}\nonumber
  \|F(\cdot;x_3,y_3)\|_{L^q(\real^2)} &\le& C 
  \|\oomega(\cdot,y_3)\|_{L^2(m)} \|b^{-m}\|_{L^2(\real^2)} 
  \|f_{x_3-y_3}\|_{L^q(\real^2)} \\ \label{Hbdd2} 
  &\le& \frac{C \|\oomega(\cdot,y_3)\|_{L^2(m)}}{|x_3-y_3|^{2-\frac2q}}~,
\end{eqnarray}
where $2-\frac2q > 1$. Combining \reff{Hbdd1}, \reff{Hbdd2}, we
obtain \reff{Heq}. \QED

\medskip
An easy consequence is the analogue of Corollary~\ref{2Dprod}: 

\begin{corollary}\label{3Dprod}
Assume that $\oomega_1, \oomega_2 \in X^2(m)$ for some $m > 1/2$, 
and let $\uu_1$ be the velocity field obtained from $\oomega_1$
via the Biot-Savart law \reff{3DBS}. Then $\uu_1 \oomega_2 \in 
X^p(m)$ for all $p \in (1,2)$, and there exists $C(m,p) > 0$ 
such that 
\begin{equation}\label{3Duomega}
  \|\uu_1 \oomega_2\|_{X^p(m)} \le C \|\oomega_1\|_{X^2(m)}
  \|\oomega_2\|_{X^2(m)}~.
\end{equation}
\end{corollary}

\proof 
Assume that $1 < p < 2$ and $\frac1q = \frac1p -\frac12$. By H\"older's 
inequality, we have for all $x_3 \in \real$
$$
  \|b(\cdot)^{m}\uu_1(\cdot,x_3)\oomega_2(\cdot,x_3)\|_{L^p(\real^2)}  
  \,\le\, \|\uu_1(\cdot,x_3)\|_{L^q(\real^2)} 
  \|b(\cdot)^m \oomega_2(\cdot,x_3)\|_{L^2(\real^2)}~.
$$
Taking the supremum over $x_3$ and using 
Proposition~\ref{velocity-vorticity}, we obtain \reff{3Duomega}. 
Moreover, since $\uu_1 \in X^q(0)$ and $\oomega_2 \in X^2(m)$, 
it is clear that $x_3 \mapsto \uu_1(\cdot,x_3)\oomega_2(\cdot,x_3)$ 
is continuous from $\real$ into $L^p(m)$. \QED

\medskip
In the rest of this section, we fix some $\lambda \in [0,1)$. Given 
$\rho \in \real$ and $\phi \in C^1_b(\real)$, our
goal is to compare the velocity field $\UU^B(\xp;\rho+\phi(x_3))$
defined by \reff{OmBUB} with the velocity field $\tilde \UU^B(x;\rho,
\phi)$ obtained from $\OO^B(\xp;\rho+\phi(x_3))$ via the Biot-Savart 
law. (As in Section~\ref{NAS_stab} we omit the dependence on $\lambda$ 
for simplicity.) Since $\OO^B$ has only the third component nonzero, 
\reff{3DBS} implies that $\tilde \UU^B$ has only the first two
components nonzero:
\begin{equation}\label{BSeq1}
  \pmatrix{\tilde U_1^B(x;\rho,\phi) \cr \tilde U_2^B(x;\rho,\phi)} 
  \,=\, \frac{1}{4\pi} \int_{\real^3} \frac{1}{|x-y|^3}
  \pmatrix{y_2-x_2 \cr x_1-y_1}\,\Omega^B(\yp;\rho+\phi(y_3))\d\yp 
  \d y_3~.
\end{equation}
On the other hand, for any $x_3 \in \real$, $\UU^B(\xp;\rho+\phi(x_3))$ 
is obtained from $\OO^B(\xp;\rho+\phi(x_3))$ via the two-dimensional
Biot-Savart law \reff{BS2Dbis}, which can be written in the form
\begin{equation}\label{BSeq2}
  \pmatrix{U_1^B(\xp;\rho+\phi(x_3)) \cr U_2^B(\xp;\rho+\phi(x_3))} 
  \,=\, \frac{1}{4\pi} \int_{\real^3} \frac{1}{|x-y|^3}
  \pmatrix{y_2-x_2 \cr x_1-y_1}\,\Omega^B(\yp;\rho+\phi(x_3))\d\yp 
  \d y_3~,
\end{equation}
because
\begin{equation}\label{2D3Did}
  \int_\real \frac{1}{|x-y|^3}\,\d y_3 \,\equiv\, \int_\real 
  \frac{1}{(|\xp-\yp|^2+(x_3-y_3)^2)^{3/2}}\d y_3 \,=\, 
  \frac{2}{|\xp-\yp|^2}~.
\end{equation}

Using these representation formulas, it is easy to show that the
velocity fields $\tilde \UU^B$, $\UU^B$ are close if the function
$\phi$ varies sufficiently slowly.

\begin{proposition}\label{UUcompare}
Fix $\lambda \in [0,1)$, and assume that $\rho \in \real$ and 
$\phi \in C^1_b(\real)$ satisfy $|\rho| + \|\phi\|_{L^\infty} 
\le R_1(\lambda)$, where $R_1$ is defined in 
Proposition~\ref{fixedpoint}. Then $\tilde \UU^B(\cdot;\rho,\phi) 
\in C^0_b(\real^3)$, and there exists $C(\lambda) > 0$ such that
\begin{eqnarray}\label{UUfirst}
  &&\sup_{x\in\real^3} |\tilde \UU^B(x;\rho,\phi)| \,\le\, 
  C(|\rho|+\|\phi\|_{L^\infty})~, \\ \label{UUcomp}
  &&\sup_{x\in\real^3} |\tilde \UU^B(x;\rho,\phi) - 
  \UU^B(\xp;\rho+\phi(x_3))| \,\le\, C \|\phi'\|_{L^\infty}~.
\end{eqnarray}
\end{proposition}

\proof Since $\Omega^B(\xp;\rho)$ is a continuous function of $\xp \in
\real^2$ which decays rapidly as $|\xp| \to \infty$, uniformly 
in $\rho \in [-R_1,R_1]$, it is not difficult to verify that
the velocity field $\tilde \UU^B(x;\rho,\phi)$ defined by \reff{BSeq1} 
depends continuously on $x \in \real^3$. Next using \reff{2D3Did} and
\reff{ombdd2} with $m = 1$, we find
$$
  |\tilde \UU^B(x;\rho,\phi)| \,\le\, \frac{1}{2\pi}\int_{\real^2}
  \frac{1}{|\xp-\yp|}\frac{C(|\rho|+\|\phi\|_{L^\infty})}{b(\yp)}
  \d\yp~.  
$$
Since $b^{-1} \in L^p(\real^2)$ for all $p \in (1,\infty]$, the
above integral is uniformly bounded for all $\xp \in \real^2$
(see Proposition~\ref{BS2Destimates}), and we obtain \reff{UUfirst}. 

Finally, taking the difference of \reff{BSeq1} and \reff{BSeq2}, we 
see that $|\tilde \UU^B - \UU^B| \le C(D_1 + D_2)$, where
$$
  D_i(x) \,=\, \int_{\real^3} \frac{|x_i-y_i|}{|x-y|^3}
  \Big|\Omega^B(\yp;\rho+\phi(y_3)) - \Omega^B(\yp;\rho+\phi(x_3))
  \Big|\d\yp \d y_3~, \quad i = 1,2~.
$$
But
\begin{eqnarray*}
  \Big|\Omega^B(\yp;\rho+\phi(y_3)) - \Omega^B(\yp;\rho+\phi(x_3))
  \Big| &\le& \Big| \int_{x_3}^{y_3} \partial_\rho 
  \Omega^B(\yp;\rho+\phi(z)) \phi'(z)\d z\Big| \\
  &\le& |x_3-y_3| \|\phi'\|_{L^\infty} \sup_{|\rho|\le R_1}
  |\partial_\rho \Omega^B(\yp;\rho)|~.
\end{eqnarray*}
Since
$$
  \int_\real \frac{|x_3-y_3|}{|x-y|^3}\,\d y_3 \,=\,\frac{2}{|\xp-\yp|}~,
  \quad \hbox{and}\quad \frac{|x_i-y_i|}{|\xp-\yp|} \,\le\, 1~,
$$
we thus find
\begin{equation}\label{Dibdd}
  \|D_i\|_{L^\infty} \,\le\, 2\|\phi'\|_{L^\infty} \int_{\real^2}
  \sup_{|\rho|\le R_1} |\partial_\rho \Omega^B(\yp;\rho)|\d\yp~,
  \quad i = 1,2~.
\end{equation}
Using now \reff{ombdd2} with $m = 2$, we see that the integrand 
is bounded by $C/b(\yp)^2$, hence the integral in \reff{Dibdd} 
is finite. This gives \reff{UUcomp}. \QED

\bigskip\noindent
{\bf Acknowledgements.}
A part of this work was completed when CEW was a visitor at Institut
Fourier, Universit\'e de Grenoble I, whose hospitality is gratefully
acknowledged. The research of CEW is supported in part by the NSF
through grant DMS-0405724, and the work of ThG is supported by the 
ACI ``Structure and dynamics of nonlinear waves'' of the French 
Ministry of Research. 


\bibliographystyle{plain}

\begin{thebibliography}{10}

\bibitem{burgers:1948}
J.~M. Burgers.
\newblock A mathematical model illustrating the theory of turbulence.
\newblock {\em Adv. Appl. Mech.}, 1:171--199, 1948.

\bibitem{crowdy:1998}
D.~G. Crowdy.
\newblock A note on the linear stability of {B}urgers vortex.
\newblock {\em Stud. Appl. Math.}, 100(2):107--126, 1998.

\bibitem{gallay:2002a}
Th. Gallay and C.~E. Wayne.
\newblock Invariant manifolds and the long-time asymptotics of the
  {N}avier-{S}tokes and vorticity equations on {$\mathbb R\sp 2$}.
\newblock {\em Arch. Ration. Mech. Anal.}, 163(3):209--258, 2002.

\bibitem{gallay:2005}
Th. Gallay and C.~E. Wayne.
\newblock Existence and stability of asymmetric {B}urgers vortices.
\newblock Preprint, 2005.

\bibitem{gallay:2004}
Th. Gallay and C.~E. Wayne.
\newblock Global stability of vortex solutions of the two-dimensional
  {N}avier-{S}tokes equation.
\newblock {\em Comm. Math. Phys.}, 255(1):97--129, 2005.

\bibitem{jimenez:1996}
J.~Jim\'enez, H.~K. Moffatt, and C.~Vasco.
\newblock The structure of the vortices in freely decaying two-dimensional
  turbulence.
\newblock {\em J. Fluid Mech.}, 313:209--222, 1996.

\bibitem{leibovich:1981}
S.~Leibovich and Ph. Holmes.
\newblock Global stability of the {B}urgers vortex.
\newblock {\em Phys. Fluids}, 24(3):548--549, 1981.

\bibitem{lieb:1997}
E.~H. Lieb and M.~Loss.
\newblock {\em Analysis}, volume~14 of {\em Graduate Studies in Mathematics}.
\newblock American Mathematical Society, Providence, RI, 1997.

\bibitem{moffatt:1994}
H.~K. Moffatt, S.~Kida, and K.~Ohkitani.
\newblock Stretched vortices---the sinews of turbulence;
  large-{R}eynolds-number asymptotics.
\newblock {\em J. Fluid Mech.}, 259:241--264, 1994.

\bibitem{prochazka:1995}
A.~Prochazka and D.~I. Pullin.
\newblock On the two-dimensional stability of the axisymmetric {B}urgers
  vortex.
\newblock {\em Phys. Fluids}, 7(7):1788--1790, 1995.

\bibitem{prochazka:1998}
A.~Prochazka and D.~I. Pullin.
\newblock Structure and stability of non-symmetric {B}urgers vortices.
\newblock {\em J. Fluid Mech.}, 363:199--228, 1998.

\bibitem{reed:1972}
M.~Reed and B.~Simon.
\newblock {\em Methods of modern mathematical physics. {I}. {F}unctional
  analysis}.
\newblock Academic Press, New York, 1972.

\bibitem{robinson:1984}
A.~C. Robinson and P.~G. Saffman.
\newblock Stability and structure of stretched vortices.
\newblock {\em Stud. Appl. Math.}, 70(2):163--181, 1984.

\bibitem{rossi:1997}
M.~Rossi and S.~Le~Diz\`es.
\newblock Three-dimensional temporal spectrum of stretched vortices.
\newblock {\em Phys. Rev. Lett.}, 78:2567--2569, 1997.

\bibitem{rudin:1991}
W.~Rudin.
\newblock {\em Functional analysis}.
\newblock International Series in Pure and Applied Mathematics. McGraw-Hill
  Inc., New York, second edition, 1991.

\bibitem{schmid:2004}
P.~J. Schmid and M.~Rossi.
\newblock Three-dimensional stability of a {B}urgers vortex.
\newblock {\em J. Fluid Mech.}, 500:103--112, 2004.

\bibitem{townsend:1951}
A.~A. Townsend.
\newblock On the fine-scale structure of turbulence.
\newblock {\em Proc. R. Soc. Lond. A}, 208:534--542, 1951.

\end{thebibliography}

\end{document}